\documentclass[10pt,a4paper]{amsart}
\usepackage{amsmath,amssymb,amsthm,enumitem}
\usepackage[a4paper]{geometry}
\usepackage{palatino}
\usepackage{graphicx}
\usepackage{float}
\usepackage{xcolor}
\usepackage[colorlinks=true]{hyperref}
\hypersetup{urlcolor=blue, citecolor=red, linkcolor=blue}
\usepackage{cite}
\IfFileExists{ulem.sty}{\usepackage[normalem]{ulem}}{\newcommand{\sout}[1]{[eliminare: ##1]}}
\newcommand{\del}[1]{{\color{red}\sout{#1}}}     

\numberwithin{equation}{section}

\newcommand{\e}{\varepsilon}
\newcommand{\R}{\mathbb{R}}
\newcommand{\dd}{\,\mathrm{d}}
\newcommand{\bd}{\boldsymbol}

\DeclareMathOperator{\sign}{sign}

\newtheorem{theorem}{Theorem}[section]
\newtheorem{lemma}[theorem]{Lemma}

\newtheorem{proposition}[theorem]{Proposition}

\theoremstyle{definition}
\newtheorem{remark}[theorem]{Remark}

\begin{document}

\title[Polygonal bubbling in an almost critical problem]{A sharp dimensional threshold for polygonal bubbling in an almost critical problem}

\author[G. M. Rago]{Giuseppe Mario Rago}
\address{\noindent  Dipartimento di Matematica, Universit\`a degli Studi di Bari Aldo Moro, Italy}
\email{g.rago6@phd.uniba.it}

\author[G. Vaira]{Giusi Vaira}
\address{\noindent  Dipartimento di Matematica, Universit\`a degli Studi di Bari Aldo Moro, Italy}
\email{giusi.vaira@uniba.it}

\author[Y. Wu]{Yue Wu}
\address{\noindent  School of Mathematical Sciences, Zhejiang Normal University, People's Republic of China}
\email{yuewu1@zjnu.edu.cn}

\thanks{G.~M. Rago and G.~Vaira acknowledge support from the INdAM-GNAMPA project ``Singolarità, nonlinearità e interazioni non locali: nuovi approcci analitici per equazioni differenziali e sistemi'' (CUP E53C25002010001).
Y.~Wu acknowledges support from the China Scholarship Council (No.~202607260015).}


\begin{abstract}
In this paper we study sign-changing solutions for the problem $-\Delta u=|u|^{2^*-2-\varepsilon}u$ in $\mathcal B$, $u=0$ on $\partial\mathcal B$, where $\mathcal B$ is the unit ball of $\mathbb R^N$, $N\geq3$, $2^*=\frac{2N}{N-2}$ and $\varepsilon>0$ is a small parameter, so that the problem is slightly subcritical. In \cite{BDP}, Bartsch, D'Aprile and Pistoia proved that, in every dimension, there exist two solutions which blow up, as $\varepsilon\to0$, positively at the origin and negatively at two antipodal points; as a by-product of our approach, we prove their conjecture that the corresponding reduced energy has exactly two critical points, both nondegenerate, in every dimension. We consider the natural generalization in which the negative bubbles are located at the vertices of a regular $k$-gon, $k\geq3$, and we show that the picture changes drastically. When $k=3$ there is a sharp dimensional threshold: if $N\geq13$ the reduced energy has exactly two critical points, a nondegenerate saddle point and a nondegenerate local minimum, which generate two distinct solutions, while if $3\leq N\leq12$ no solution of this type exists. When $k\geq4$ no solution of this type exists in any dimension. Moreover, for both families of solutions found for $k=3$ the nodal set stays uniformly away from $\partial\mathcal B$, in contrast with the case of two negative bubbles, where the nodal set of one of the two solutions touches the boundary, and with the Brezis--Nirenberg problem, where the nodal set of the solutions with the same configuration always reaches the boundary.
\end{abstract}

\maketitle

\section{Introduction}
In this paper we are interested in sign-changing solutions of the slightly subcritical problem
\begin{equation}\label{pb}
\begin{cases}
-\Delta u = |u|^{2^{*}-2-\e} u & \mbox{in } \mathcal B, \\
u=0 & \mbox{on } \partial \mathcal B,
\end{cases}
\end{equation}
where $\mathcal B := \{ x \in \R^{N} : |x| < 1\}$, $N \geq 3$, $2^{*}:= \frac{2N}{N-2}$ is the critical exponent for the Sobolev embedding and $\e$ is a small positive parameter.

The critical exponent plays a special role. If $\e=0$, an application of Poho\v{z}aev's identity \cite{Po} shows that the problem has no nontrivial solutions in $\mathcal B$, and more generally in any star-shaped domain. If instead $\e>0$, the compactness of the embedding $H^1_0(\mathcal B)\hookrightarrow L^{2^*-\e}(\mathcal B)$ yields a positive solution and infinitely many sign-changing solutions. Therefore, as $\e\to0$, the solutions of \eqref{pb} must blow up, and a natural problem is to describe their asymptotic profile.

For positive solutions on a general smooth bounded domain $\Omega$ the picture is rather complete. The least energy solution blows up at a single point, which is a critical point of the Robin function of $\Omega$, and its profile is a standard bubble, see \cite{BP,H,R,FW}. Positive solutions blowing up at several points have been studied in \cite{BLR,R2}; in particular, they do not exist if $\Omega$ is convex, see \cite{GrTa}, while on an annulus multi-peak positive solutions have been constructed in \cite{MRV}, by exploiting the explicit form of the Green and Robin functions.

Much less is known for sign-changing solutions. Low energy nodal solutions were classified in \cite{BEP}, according to the concentration speeds of their positive and negative parts. Solutions with one positive and one negative bubble concentrating at two different points were built in \cite{BMP}, where the location of their nodal set is also investigated, and towers of bubbles of alternating sign concentrating at a single point were constructed in \cite{PW,MuPi}; four-bubble nodal solutions in convex symmetric domains were obtained in \cite{BDP2}. In all these constructions the main difficulty is not the Lyapunov--Schmidt reduction itself, which is by now standard, but the analysis of the critical points of the resulting finite-dimensional function, in which the location of the bubbles and their concentration rates are coupled.

Our starting point is the paper by Bartsch, D'Aprile and Pistoia \cite{BDP}. They consider problem \eqref{pb} in the ball and, for every $N\geq3$, they find two distinct solutions which, as $\e\to0$, look like
\[
P\mathcal U_{\delta_0,0}-P\mathcal U_{\delta_1,\xi}-P\mathcal U_{\delta_1,-\xi},
\]
that is, a positive bubble at the origin and two negative bubbles at two antipodal points $\pm\xi$ (we refer to Section~\ref{sec:main} for the notation). The two solutions correspond to a local minimum and a local saddle point of the reduced energy, so that they are not detected by a global min--max scheme, and they are distinguished by the distance of the negative peaks from the origin. Remarkably, they also have different nodal structures: the nodal set of the solution whose negative peaks are farther from the origin touches $\partial\mathcal B$, whereas the nodal set of the other one is a closed surface around the positive bubble, which does not touch $\partial\mathcal B$. In \cite[Remark~3.2]{BDP} it is conjectured that the reduced energy has exactly these two critical points, and their nondegeneracy is checked numerically for $N<6$; as a by-product of our approach we prove this conjecture in every dimension $N\geq3$ (see Remark~\ref{rmk:k2}).

The same kind of configuration has recently been studied for the Brezis--Nirenberg problem
\begin{equation}\label{BNpb}
-\Delta u=|u|^{2^*-2}u+\e u\quad \mbox{in } \mathcal B,\qquad u=0\quad \mbox{on } \partial\mathcal B,
\end{equation}
which is very close to \eqref{pb}: in both problems the solutions blow up as $\e\to0$, the profile near each blow-up point is a bubble, and the location of the blow-up points is governed by the Green and the Robin functions of the domain. In dimension $N=4$, solutions of \eqref{BNpb} which blow up positively at the centre of the ball and negatively at the $k$ vertices of a regular polygon were constructed in \cite{PRV2}, and their persistence in higher dimensions has been studied in \cite{PRV3}. The answer depends on the dimension: if $N=5$ such solutions exist for every $k\geq2$, and the ring of negative peaks approaches $\partial\mathcal B$ as $k\to+\infty$; if $N\geq6$ there are two solutions for $2\leq k\leq6$, while no solution of this type exists for $k$ large. Moreover, for $2\leq k\leq6$ the nodal set of every such solution reaches $\partial\mathcal B$. In particular, already for $k=2$ the Brezis--Nirenberg problem behaves differently from \eqref{pb}: none of its solutions has a nodal set consisting of a closed surface around the positive bubble, as happens for one of the solutions in \cite{BDP}.
In view of these results, it is natural to ask what happens for the almost critical problem \eqref{pb} when the two antipodal points of \cite{BDP} are replaced by the $k\geq3$ vertices of a regular polygon. More precisely, we look for solutions of \eqref{pb} which, as $\e \to 0$, look like a positive bubble located at the origin and $k \geq 3$ negative bubbles located at the vertices $\xi_1(\rho),\dots,\xi_k(\rho)$ of a regular $k$-gon of radius $\rho\in(0,1)$, namely
\begin{equation}\label{ansatz-intro}
u_\e=P\mathcal U_{\delta_0,0}-\sum_{j=1}^kP\mathcal U_{\delta_1,\xi_j(\rho)}+\phi_\e,
\end{equation}
where $\phi_\e$ is a small remainder term. Our aim is to understand for which pairs $(N,k)$ these configurations occur and to describe the behaviour of their nodal sets near the boundary.

Given the analogy between \eqref{pb} and \eqref{BNpb}, one might expect a picture similar to the one found in \cite{PRV3}. Our results show that this is not the case, and that the almost critical problem is much more rigid.

To state them, we introduce some notation, referring to Section~\ref{sec:main} for the details. For $\delta>0$ and $\xi\in\R^N$, $\mathcal U_{\delta,\xi}$ is the standard bubble defined in \eqref{bubble} and $P\mathcal U_{\delta,\xi}$ is its projection onto $H^1_0(\mathcal B)$. For $\rho\in(0,1)$, $\xi_1(\rho),\dots,\xi_k(\rho)$ are the vertices of the regular $k$-gon of radius $\rho$ in the plane $\{x_3=\dots=x_N=0\}$ defined in \eqref{xij}, and $\mathcal H_k$ is the space of the functions of $H^1_0(\mathcal B)$ which have the symmetries of this polygon in the $(x_1,x_2)$-plane and are radial with respect to $(x_3,\dots,x_N)$, see \eqref{Hk}. Finally, the concentration parameters in \eqref{ansatz-intro} are written as
\begin{equation}\label{pallino}
\delta_0=(\mathtt d_N\,\e)^{\frac1{N-2}}\,d_0^{\frac2{N-2}},\qquad \delta_1=(\mathtt d_N\,\e)^{\frac1{N-2}}\,d_1^{\frac2{N-2}},\qquad (d_0,d_1,\rho)\in\mathcal M:=(0,+\infty)^2\times(0,1),
\end{equation}
where $\mathtt d_N>0$ is the constant defined in \eqref{dN}, and we set $\mathcal V_\e:=P\mathcal U_{\delta_0,0}-\sum_{j=1}^kP\mathcal U_{\delta_1,\xi_j(\rho)}$.

Our first result gives a complete picture in the case of three negative bubbles, with the sharp dimensional threshold $N=13$.

\begin{figure}[t]
\centering
\includegraphics[width=0.4\textwidth]{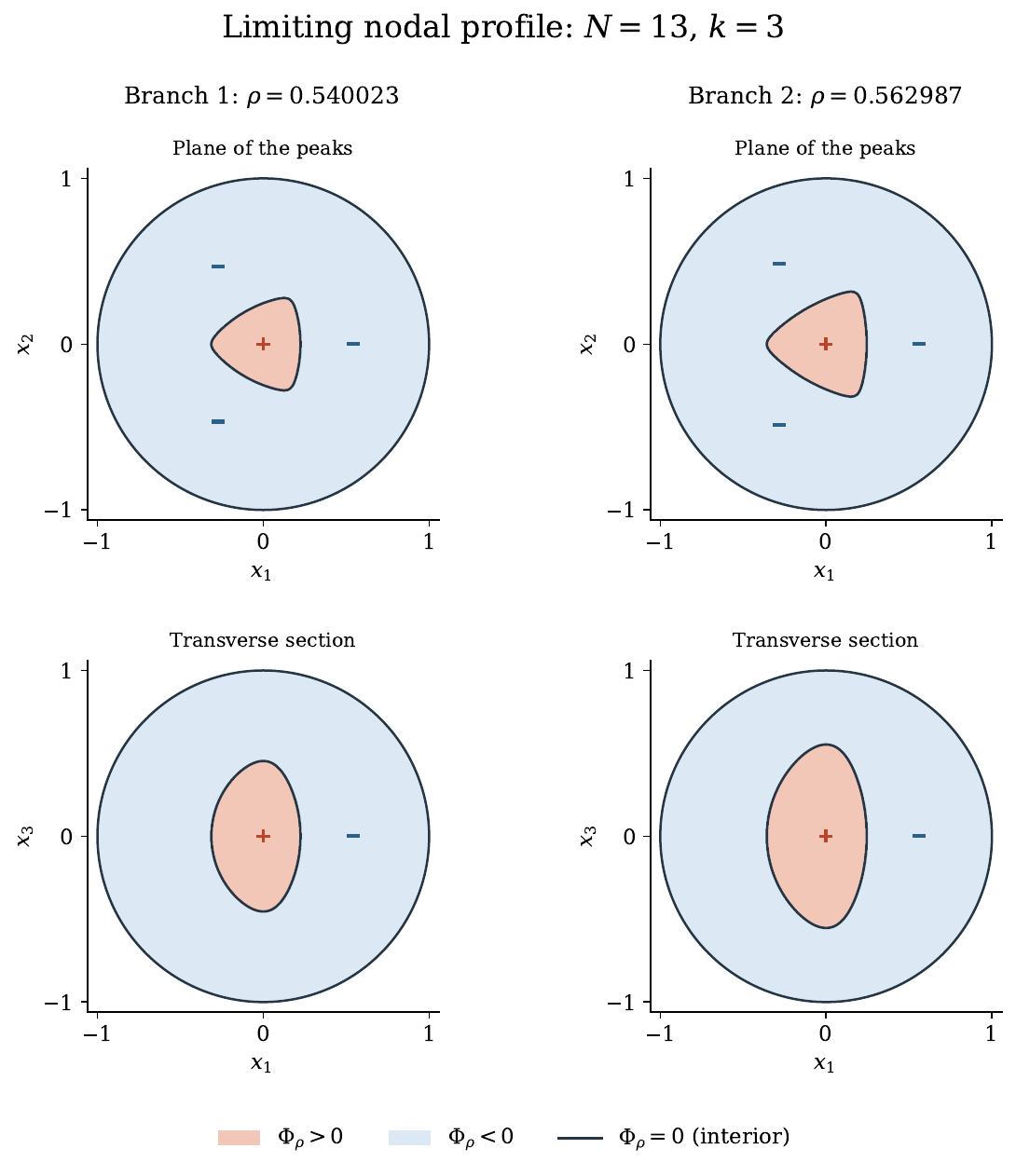}
\caption{\tiny{Sign of the limit profile $\Phi_\rho$ in \eqref{eq:nodal-profile} for $N=13$, $k=3$, at the critical radii $\rho^*\approx0.5400$ (left) and $\rho^{**}\approx0.5630$ (right). Top: the plane of the peaks $\{x_3=\dots=x_N=0\}$; bottom: the transverse section $\{x_2=x_4=\dots=x_N=0\}$. The marks $+$ and $-$ denote the positive and the negative concentration points. In both cases the nodal set $\{\Phi_\rho=0\}$ is a closed hypersurface around the origin, uniformly separated from $\partial\mathcal B$.}}
\label{fig:nodal}
\end{figure}

\begin{theorem}\label{thm:main}
Let $k=3$.
\begin{enumerate}[label=(\roman*)]
\item\emph{(\textbf{Existence})} If $N\geq13$, then the reduced function $\Psi_3$ defined in \eqref{psik} has exactly two critical points in $\mathcal M$, both nondegenerate: a saddle point of Morse index $1$ and a strict local minimum. Consequently, there exists $\e_0>0$ such that for every $\e\in(0,\e_0)$ problem \eqref{pb} has two distinct sign-changing solutions $u^{(1)}_\e,u^{(2)}_\e\in \mathcal H_3$ of the form
\[
u^{(i)}_\e=P\mathcal U_{\delta^{(i)}_{0,\e},0}-\sum_{j=1}^3P\mathcal U_{\delta^{(i)}_{1,\e},\xi_j(\rho^{(i)}_\e)}+\phi^{(i)}_\e,\qquad i=1,2,
\]
where $\delta^{(i)}_{0,\e}$ and $\delta^{(i)}_{1,\e}$ are given by \eqref{pallino} with $(d_0,d_1)=(d^{(i)}_{0,\e},d^{(i)}_{1,\e})$, $\|\phi^{(i)}_\e\|\to0$ and, as $\e\to0$,
\[
\big(d^{(1)}_{0,\e},d^{(1)}_{1,\e},\rho^{(1)}_\e\big)\to\big(\bd d_0(\rho^*),\bd d_1(\rho^*),\rho^*\big),\qquad
\big(d^{(2)}_{0,\e},d^{(2)}_{1,\e},\rho^{(2)}_\e\big)\to\big(\bd d_0(\rho^{**}),\bd d_1(\rho^{**}),\rho^{**}\big).
\]
Here $\bd d_0(\rho)$ and $\bd d_1(\rho)$ are the explicit functions given in \eqref{solsis} (with $k=3$), $\rho_0:=\rho_{0,3}$ is the unique zero of the function $A_3$ defined in \eqref{Akrho}, see \eqref{imprho0}, and the critical radii $\rho^*,\rho^{**}$ satisfy
\[
\rho_0<\rho^*<\frac{11}{20}<\rho^{**}<\frac1{\sqrt3}.
\]
The two branches of solutions and their nodal regions are illustrated, for $N=13$, in Figure~\ref{fig:nodal}.
\item\emph{(\textbf{Nonexistence})} If $3\leq N\leq12$, then $\Psi_3$ has no critical points in $\mathcal M$. Consequently, there is no family of solutions $u_\e\in \mathcal H_3$ of \eqref{pb} of the form
\[
u_\e=\mathcal V_\e+\phi_\e,\qquad \phi_\e\in\mathcal K^\perp_{d_{0,\e},d_{1,\e},\rho_\e},\quad \|\phi_\e\|\to0,
\]
where $\mathcal K^\perp_{d_0,d_1,\rho}$ is the space of functions satisfying the orthogonality conditions of the Lyapunov--Schmidt reduction, see Proposition~\ref{prop:reduction},
with $(d_{0,\e},d_{1,\e},\rho_\e)\to p\in\mathcal M$ as $\e\to0$.
\end{enumerate}
\end{theorem}

For four or more negative bubbles, instead, configurations of this type never occur, in any dimension.

\begin{theorem}\label{thm:k-geq4}
Let $N\geq3$ and $k\geq4$. Then the reduced function $\Psi_k$ defined in \eqref{psik} has no critical points in $\mathcal M$. Consequently, there is no family of solutions $u_\e\in\mathcal H_k$ of \eqref{pb} of the form
\[
u_\e=\mathcal V_\e+\phi_\e,\qquad
\phi_\e\in\mathcal K^\perp_{d_{0,\e},d_{1,\e},\rho_\e},\quad
\|\phi_\e\|\to0,
\]
with $(d_{0,\e},d_{1,\e},\rho_\e)\to p\in\mathcal M$ as $\e\to0$.
\end{theorem}
Finally, the two solutions found for $k=3$ are negative in a whole neighbourhood of $\partial\mathcal B$, so that their nodal sets do not touch the boundary.

\begin{theorem}\label{thm:nodal}
Let $N\geq13$, and let $u_\e^{(1)},u_\e^{(2)}$ be the two
solution families of Theorem~\ref{thm:main}.
For a solution $u$ of \eqref{pb}, let $\mathcal Z(u):=\overline{\{x\in\mathcal B:u(x)=0\}}^{\,\overline{\mathcal B}}$ denote the closure of its internal nodal set. There exist $\tau\in(0,1)$ and $\e_1>0$ such that, for
$i=1,2$ and $0<\e<\e_1$,
\[
\partial_\nu u_\e^{(i)}>0\quad\text{on }\partial\mathcal B,
\qquad
u_\e^{(i)}(x)<0\quad\text{if }1-\tau<|x|<1,
\]
where $\nu$ is the exterior unit normal. In particular,
\[
\mathcal Z(u_\e^{(i)})\subset\{x\in\mathcal B:|x|\leq1-\tau\},
\]
so the internal nodal sets of both families remain uniformly
separated from $\partial\mathcal B$ as $\e\to0$.
\end{theorem}

Let us first compare these results with those of \cite{BDP}. As in the case $k=2$ (see Remark~\ref{rmk:k2}), when solutions exist they come in a pair, generated by a nondegenerate local minimum and a nondegenerate saddle point of the reduced energy, which has no other critical points. However, while two antipodal negative bubbles are admissible in every dimension $N\geq3$, three negative bubbles are admissible only if $N\geq13$, and four or more are never admissible. Moreover, for $k=2$ the nodal set of one of the two solutions reaches $\partial\mathcal B$, while for $k=3$ this never happens: both solutions behave like the solution of \cite{BDP} whose negative peaks are closer to the origin.

The comparison with the Brezis--Nirenberg problem is even more striking, since the two problems turn out to behave in opposite ways. First, \eqref{BNpb} admits rings of any small number of negative peaks, $2\leq k\leq6$, in every dimension $N\geq5$, and rings with any number of peaks if $N=5$, whereas \eqref{pb} admits at most three of them. Second, the dimension plays opposite roles: for \eqref{BNpb} increasing the dimension reduces the admissible values of $k$, while for \eqref{pb} a high dimension is necessary to have three negative bubbles. Third, and most importantly, the nodal sets behave in opposite ways. For $k=3$ and $N\geq13$ both problems have two solutions of the form \eqref{ansatz-intro} (for \eqref{pb} they are generated by the only two critical points of the reduced energy $\Psi_3$, and every solution of this form concentrates at one of them, see Remark~\ref{rmk:limit}); however, the nodal sets of the two solutions of \eqref{BNpb} always reach $\partial\mathcal B$, since these solutions are positive near the $(N-3)$-dimensional sphere $\partial\mathcal B\cap\{x_1=x_2=0\}$ and negative near the points $\xi_j/|\xi_j|$ (see \cite[Theorem~1.5]{PRV3}), while the two solutions of \eqref{pb} are negative in a whole neighbourhood of $\partial\mathcal B$. We point out that, in the plane of the peaks, the nodal set of one of the solutions of \cite{PRV3} is a closed curve around the origin, but it reaches $\partial\mathcal B$ outside this plane; for \eqref{pb}, instead, the nodal set is a genuinely closed hypersurface around the positive bubble. Hence the same configuration of bubbles produces, in the two problems, solutions with opposite nodal behaviour near the boundary.

This difference can be explained as follows. In both problems the sign of the exterior normal derivative on $\partial\mathcal B$ is decided, as $\e\to0$, by the comparison between the ratio of the concentration rates of the positive and of the negative bubbles and a boundary density which depends on the radius of the polygon (see \cite[Remark~1.6]{PRV3} and Proposition~\ref{prop:nodal-normal}). In \eqref{BNpb} the concentration rates are of order $\e^{\frac1{N-4}}$, their self-interaction is of the same order as the interaction between different bubbles, and the radius and the scales must be determined together; at every limiting radius, this comparison changes sign along $\partial\mathcal B$. In \eqref{pb}, instead, the concentration rates are of order $\e^{\frac1{N-2}}$ and the reduced energy has a logarithmic structure, so that the equations for the scales can be solved explicitly for every fixed radius. At a critical radius the ratio of the scales coincides with an explicit quantity which describes the interaction of each negative bubble with the boundary and with the other negative bubbles, and this forces the normal derivative to be positive on the whole $\partial\mathcal B$.

As in \cite{PRV3}, the proofs rely on a Lyapunov--Schmidt reduction, which reduces the problem to the study of the critical points of an explicit function $\Psi_k$ of the three variables $(d_0,d_1,\rho)$, where $d_0$ and $d_1$ are the scaled concentration parameters of the positive and of the negative bubbles, respectively. The dependence of $\Psi_k$ on $\rho$ is encoded in two coefficients: $\beta(\rho)$, which measures the interaction between the central bubble and each negative bubble, and $A_k(\rho)$, which collects the interaction of each negative bubble with the boundary and with the other negative bubbles. For every admissible $\rho$ we solve explicitly the equations in $d_0$ and $d_1$ and we show that the corresponding partial Hessian is positive definite; this reduces the problem to the study of a function $f_k$ of the single variable $\rho$, whose critical points have the same nondegeneracy and Morse index properties as those of $\Psi_k$. The delicate part is the study of the sign of $f_k'$, which requires estimates uniform in the dimension and, for $k=3$, sharp enough to capture the threshold $N=13$: existence for $N\geq13$ follows from explicit estimates of $f_3'$ at the radius $\frac{11}{20}$, together with a monotonicity argument based on a power series expansion, which gives the exact number of critical points; nonexistence for $3\leq N\leq12$ and for $k\geq4$ is obtained by showing that $f_k$ is strictly increasing. Finally, the information on the nodal sets is derived from the limit profile of the solutions away from the concentration points, which is a signed combination of Green's functions.

The paper is organized as follows. In Section~\ref{sec:main} we introduce the setting and the reduced function, and we recall the Lyapunov--Schmidt reduction. In Section~\ref{sec:onedim} we reduce the study of $\Psi_k$ to a one-dimensional problem. Section~\ref{sec:k3} is devoted to the case $k=3$ and contains the proof of Theorem~\ref{thm:main}. In Section~\ref{nodal} we study the nodal sets and prove Theorem~\ref{thm:nodal}, and in Section~\ref{nE4} we prove Theorem~\ref{thm:k-geq4}. Some technical lemmas are collected in the Appendix.

%
\section{Setting and the reduced problem}\label{sec:main}
\subsection{Setting}
Throughout the paper $G(x,y)$ denotes the Green's function of $-\Delta$ in $\mathcal B$ with Dirichlet boundary conditions, normalized so that
\begin{equation}\label{Green}
G(x,y)=\frac{1}{|x-y|^{N-2}}-H(x,y),\qquad
H(x,y)=\frac{1}{\big(|x|^2|y|^2-2\langle x,y\rangle+1\big)^{\frac{N-2}{2}}},
\end{equation}
where $H$ is the regular part of $G$.

For $\delta>0$ and $\xi\in\R^N$ we consider the standard bubbles
\begin{equation}\label{bubble}
\mathcal U_{\delta,\xi}(x):=\alpha_N\frac{\delta^{\frac{N-2}{2}}}{\big(\delta^2+|x-\xi|^2\big)^{\frac{N-2}{2}}},\qquad \alpha_N:=\big(N(N-2)\big)^{\frac{N-2}{4}},
\end{equation}
which are all (see \cite{A, CGS, T}) the positive solutions of
$$-\Delta \mathcal U=\mathcal U^{2^*-1} \ \ \mbox{in} \ \ \R^N$$
 and we denote by $P\mathcal U_{\delta,\xi}$ the projection of $\mathcal U_{\delta,\xi}$ onto $H^1_0(\mathcal B)$, that is
 $$\Delta P\mathcal U_{\delta,\xi}=\Delta \mathcal U_{\delta,\xi} \ \ \mbox{in} \ \ \mathcal B \ \ \mbox{and} \ \ \ P\mathcal U_{\delta,\xi}=0 \ \ \mbox{on} \ \ \partial\mathcal B.$$

\noindent
Let us consider an integer $k\geq3$ and let us write $x=(x',x'')\in\R^2\times\R^{N-2}$. For every  $\rho\in(0,1)$ we place the vertices of a regular $k$-gon of radius $\rho$ in the $x'$-plane, namely
\begin{equation}\label{xij}
\xi_j(\rho)=\xi_{j,k}(\rho):=\rho\Big(\cos\frac{2\pi(j-1)}{k},\ \sin\frac{2\pi(j-1)}{k},\ 0,\dots,0\Big)\in \R^2\times\R^{N-2},\qquad j=1,\dots,k.
\end{equation}
We now introduce the symmetric functional setting. We consider the orthogonal transformations of $\R^N$
\[
\mathcal R_k(x',x''):=\Big(x_1\cos\frac{2\pi}{k}-x_2\sin\frac{2\pi}{k},\ x_1\sin\frac{2\pi}{k}+x_2\cos\frac{2\pi}{k},\ x''\Big),
\qquad
\mathcal S(x',x''):=(x_1,-x_2,x''),
\]
and, for every $\mathcal O\in O(N-2)$,
\[
\widehat{\mathcal O}(x',x''):=(x',\mathcal O x'').
\]
Let $\Gamma_k$ be the subgroup of $O(N)$ generated by $\mathcal R_k$, $\mathcal S$ and $\{\widehat{\mathcal O}:\ \mathcal O\in O(N-2)\}$, and set
\begin{equation}\label{Hk}
\mathcal H_k:=\Big\{u\in H^1_0(\mathcal B):\ u\circ\gamma=u\ \text{ a.e. in }\mathcal B\ \text{ for every }\gamma\in\Gamma_k\Big\}.
\end{equation}
Note that $\Gamma_k=D_k\times O(N-2)$, where $D_k$ is the dihedral group of order $2k$ acting on $x'$, i.e. the symmetry group of the regular $k$-gon with vertices $\xi_1(1),\dots,\xi_k(1)$ defined in \eqref{xij}. Thus $\mathcal H_k$ consists of the functions which have the symmetries of this polygon in the $x'$-plane and which are radial with respect to $x''$.
Since every $\gamma\in\Gamma_k$ is an orthogonal map with $\gamma(\mathcal B)=\mathcal B$, the map $u\mapsto u\circ\gamma$ is a linear isometry of $H^1_0(\mathcal B)$; hence $\mathcal H_k$ is a closed subspace of $H^1_0(\mathcal B)$, and it is a Hilbert space with the scalar product $\int_{\mathcal B}\nabla u\cdot\nabla v\dd x$. Moreover, the energy functional
\[
J_\e(u):=\frac12\int_{\mathcal B}|\nabla u|^2\dd x-\frac{1}{2^*-\e}\int_{\mathcal B}|u|^{2^*-\e}\dd x
\]
satisfies $J_\e(u\circ\gamma)=J_\e(u)$ for every $\gamma\in\Gamma_k$, so that, by the principle of symmetric criticality, critical points of $J_\e$ restricted to $\mathcal H_k$ are solutions of \eqref{pb}.

We look for solutions of the form
\begin{equation}\label{ansatz}
u_\e=\mathcal V_\e+\phi_\e,\qquad
\mathcal V_\e:=P\mathcal U_{\delta_0,0}-\sum_{j=1}^kP\mathcal U_{\delta_1,\xi_j(\rho)},
\end{equation}
where the concentration parameters are
\begin{equation}\label{deltas}
\delta_0:=(\mathtt d_N\,\e)^{\frac1{N-2}}\,d_0^{\frac2{N-2}},\qquad \delta_1:=(\mathtt d_N\,\e)^{\frac1{N-2}}\,d_1^{\frac2{N-2}},\qquad (d_0,d_1,\rho)\in\mathcal M:=(0,+\infty)^2\times(0,1),
\end{equation}
$\phi_\e$ is a small remainder and $\mathtt d_N>0$ is the constant
\begin{equation}\label{dN}
\mathtt d_N:=\frac{\omega_N}{\alpha_N\displaystyle\int_{\R^N}\mathcal U_{1,0}^{2^*-1}\dd x},\qquad \omega_N:=\frac1{2^*}\int_{\R^N}\mathcal U_{1,0}^{2^*}\dd x.
\end{equation}
This choice of $\mathtt d_N$ makes the interaction terms and the logarithmic terms in the expansion of the energy appear with the same coefficient $\omega_N$, see \eqref{expansion}. Note that $\mathcal V_\e\in \mathcal H_k$, since the configuration $\{0\}\cup\{\xi_1(\rho),\dots,\xi_k(\rho)\}$ is invariant under the symmetries defining $\mathcal H_k$, and that all the negative bubbles have the same concentration rate $\delta_1$.

\subsection{The reduced function}
Set
\begin{equation}\label{betarho}
\beta(\rho) := G(0, \xi_{1}(\rho)) = \rho^{2-N}-1
\end{equation}
and
\begin{equation}\label{Akrho}
A_{k}(\rho) := H(\xi_{1}(\rho), \xi_{1}(\rho)) - \sum_{j=2}^{k} G(\xi_1(\rho), \xi_{j}(\rho)).
\end{equation}
The expansion of the energy of the ansatz \eqref{ansatz} (up to the multiplicative constant $\omega_N$ defined in \eqref{dN}, see Proposition~\ref{prop:reduction} below) leads to the function
\begin{equation}\label{psik1}
\begin{aligned}
\Psi_{k}(d_0, d_1, \rho) =  & \frac 12 d_0^{2}H(0,0) + \frac k2 d_1^{2} H(\xi_1, \xi_1) + k d_0 d_1 G(0,\xi_1) 
 - d_1^{2} \sum_{1 \leq i <j\leq k} G(\xi_i, \xi_j) - \log d_0 - k\log d_1.
\end{aligned}
\end{equation}
Now, $H(0,0)=1$ and, by the rotational symmetry,
$$\sum_{1\leq i<j \leq k} G(\xi_i, \xi_j) = \frac k2 \sum_{j=2}^{k} G(\xi_1, \xi_j).$$
Therefore, using the definitions of $\beta(\rho)$ and $A_{k}(\rho)$ given in \eqref{betarho} and \eqref{Akrho}, we can rewrite \eqref{psik1} as
\begin{equation}\label{psik}
\Psi_{k}(d_0, d_1, \rho) = \frac 12 d_0^{2} + \frac{k}{2} A_{k}(\rho) d_1^{2} + k\beta(\rho) d_0 d_1 - \log d_0 -k\log d_1,
\end{equation}
for $(d_0, d_1, \rho) \in \mathcal M$. We refer to $\Psi_k$ as the reduced function associated with one positive bubble and $k$ negative bubbles.

The link between $\Psi_k$ and problem \eqref{pb} is given by the following finite-dimensional reduction, which is by now standard (see e.g. \cite{MP,BMP,BDP} and the references therein). For $(d_0,d_1,\rho)\in\mathcal M$ let
\[
\mathcal K_{d_0,d_1,\rho}:=\mathrm{span}\big\{\partial_{d_0}\mathcal V_\e,\ \partial_{d_1}\mathcal V_\e,\ \partial_\rho \mathcal V_\e\big\}\subset \mathcal H_k,
\]
let $\mathcal K^\perp_{d_0,d_1,\rho}$ be its orthogonal complement in $\mathcal H_k$ with respect to the scalar product $\int_{\mathcal B}\nabla u\cdot\nabla v$, and let $\Pi^\perp$ be the corresponding orthogonal projection.

\begin{proposition}\label{prop:reduction}
For every compact set $\mathcal C\subset\mathcal M$ there exist $\e_{\mathcal C}>0$ and $\eta_{\mathcal C}>0$ such that, for every $\e\in(0,\e_{\mathcal C})$ and every $(d_0,d_1,\rho)\in\mathcal C$, there exists a unique $\phi=\phi_\e(d_0,d_1,\rho)\in\mathcal K^\perp_{d_0,d_1,\rho}$ with $\|\phi\|\leq\eta_{\mathcal C}$ such that
\[
\Pi^\perp\big[\nabla J_\e\big(\mathcal V_\e(d_0,d_1,\rho)+\phi\big)\big]=0.
\]
Moreover $\|\phi_\e\|\to0$ as $\e\to0$ uniformly on $\mathcal C$, the map $(d_0,d_1,\rho)\mapsto\phi_\e(d_0,d_1,\rho)$ is of class $C^1$, and $\mathcal V_\e+\phi_\e$ is a critical point of $J_\e$ (and hence a solution of \eqref{pb}) if and only if $(d_0,d_1,\rho)$ is a critical point of
\[
\widetilde J_\e(d_0,d_1,\rho):=J_\e\big(\mathcal V_\e+\phi_\e\big).
\]
Finally, there exist constants $\mathtt a_N,\mathtt b_N,\mathtt c_N\in\R$, depending only on $N$ and $k$, such that
\begin{equation}\label{expansion}
\widetilde J_\e(d_0,d_1,\rho)=\mathtt a_N+\mathtt b_N\,\e+\mathtt c_N\,\e\log\e+\omega_N\,\e\,\Psi_k(d_0,d_1,\rho)+o(\e)
\end{equation}
in the $C^1$-sense, uniformly with respect to $(d_0,d_1,\rho)\in\mathcal C$.
\end{proposition}

\begin{proof}[Comments on the proof]
Proposition~\ref{prop:reduction} is a particular case of the finite-dimensional reduction for the slightly subcritical problem with bubbles of prescribed signs, developed in \cite{BMP,BDP2} (see also \cite{BLR,DFM1,DFM2}), which was used in \cite{BDP} for the ball in the space of functions even in each variable. More precisely, in \cite[Lemmas~2.1 and 2.2 and Proposition~2.3]{BDP2} the following facts are proved for a smooth bounded domain, for $m$ bubbles $P\mathcal U_{\delta_i,\xi_i}$ with prescribed signs $a_i\in\{-1,1\}$, and for parameters in a compact set:
\begin{itemize}[leftmargin=2em]
\item the existence and uniqueness of a small remainder, orthogonal to the derivatives of the bubbles with respect to all the concentration parameters and to all the coordinates of the concentration points, which solves the projected equation and depends in a $C^1$ way on the parameters;
\item the equivalence between critical points of the reduced energy and solutions of the problem, which goes back to \cite{BLR};
\item the expansion of the reduced energy, uniform in the $C^1$-sense, whose principal part, when $\delta_i^{N-2}=\mathtt d_N\,\e\,\Lambda_i^2$, is $\omega_N\e$ times
\[
\frac12\sum_{i=1}^m\Lambda_i^2H(\xi_i,\xi_i)-\sum_{1\leq i<j\leq m}a_ia_j\Lambda_i\Lambda_jG(\xi_i,\xi_j)-\log(\Lambda_1\cdots\Lambda_m).
\]
\end{itemize}
We apply these results with $m=k+1$, $a_0=1$, $a_1=\dots=a_k=-1$, $\xi_0=0$, $\xi_j=\xi_j(\rho)$, $\Lambda_0=d_0$ and $\Lambda_1=\dots=\Lambda_k=d_1$; with these choices the principal part above is exactly \eqref{psik1}, that is, $\Psi_k$. Since all the operators involved commute with the action of $\Gamma_k$, the proof of \cite[Lemma~2.1]{BDP2} can be carried out in $\mathcal H_k$. In this space the derivatives of the bubbles which are not $\Gamma_k$-invariant, namely those with respect to the coordinates of the central point and those in the tangential and $x''$-directions at the vertices, disappear, and the $\Gamma_k$-invariant part of the span of the remaining ones is exactly $\mathcal K_{d_0,d_1,\rho}$. This gives the first part of the statement, and the argument of \cite{BLR} gives the equivalence between critical points of $\widetilde J_\e$ and solutions of \eqref{pb}. Moreover, by uniqueness, for these symmetric configurations the remainder provided by \cite[Lemma~2.1]{BDP2} is $\Gamma_k$-invariant and coincides with $\phi_\e(d_0,d_1,\rho)$. Hence $\widetilde J_\e$ is the composition of the reduced energy of \cite{BDP2} with the smooth map $(d_0,d_1,\rho)\mapsto(\Lambda,\xi)$ defined above, and \eqref{expansion} follows from \cite[Proposition~2.3]{BDP2} and the chain rule.
\end{proof}

By Proposition~\ref{prop:reduction}, the construction of solutions of the form \eqref{ansatz} reduces to finding stable critical points of $\Psi_k$ in $\mathcal M$. Conversely, every family of solutions of the form \eqref{ansatz} whose parameters converge in $\mathcal M$ concentrates at a critical point of $\Psi_k$.

\begin{lemma}\label{lem:limit-config}
Let $N\geq3$ and $k\geq2$. Let $\e_n\to0$ and let $u_n=\mathcal V_{\e_n}(p_n)+\phi_n\in\mathcal H_k$ be solutions of \eqref{pb} with $\e=\e_n$, where $p_n=(d_{0,n},d_{1,n},\rho_n)\in\mathcal M$, $\phi_n\in\mathcal K^\perp_{p_n}$ and $\|\phi_n\|\to0$. If $p_n\to p\in\mathcal M$, then $\nabla\Psi_k(p)=0$.
\end{lemma}
\begin{proof}
Let $\mathcal C\subset\mathcal M$ be a compact neighbourhood of $p$. For $n$ large we have $p_n\in\mathcal C$, $\e_n<\e_{\mathcal C}$ and $\|\phi_n\|\leq\eta_{\mathcal C}$. Since $\nabla J_{\e_n}(u_n)=0$, in particular $\Pi^\perp\nabla J_{\e_n}(\mathcal V_{\e_n}(p_n)+\phi_n)=0$, and the uniqueness assertion in Proposition~\ref{prop:reduction} gives $\phi_n=\phi_{\e_n}(p_n)$. Hence, again by Proposition~\ref{prop:reduction}, $p_n$ is a critical point of $\widetilde J_{\e_n}$, and \eqref{expansion} yields
\[
0=(\omega_N\e_n)^{-1}\nabla\widetilde J_{\e_n}(p_n)=\nabla\Psi_k(p_n)+o(1).
\]
Passing to the limit, we get $\nabla\Psi_k(p)=0$.
\end{proof}

\section{Reduction to a one-dimensional problem}\label{sec:onedim}

\subsection{Preliminary properties}
We first collect some preliminary properties of the functions $A_{k}$ and $\beta$ defined in \eqref{Akrho} and \eqref{betarho}, respectively. A direct computation shows that
$$\frac{d}{d\rho} H(\xi_1(\rho), \xi_1(\rho))>0 \ \ \ \ \mbox{and} \ \ \ \ \frac{d}{d\rho} G(\xi_1(\rho), \xi_j(\rho))<0 \ \ \ \mbox{for all} \ \ j \geq 2.$$
Therefore, by \eqref{Akrho}, $A'_{k}(\rho)>0$ for every $\rho\in(0,1)$. Recalling the expression of $G$ in \eqref{Green}, we can explicitly write
\begin{equation}\label{Akk}
A_{k}(\rho) = \frac{1}{(1-\rho^{2})^{N-2}} - \frac{1}{\rho^{N-2}} \sum_{j=1}^{k-1} \frac{1}{\left(2\sin\left(\frac{\pi j}{k} \right) \right)^{N-2}} + \sum_{j=1}^{k-1} \frac{1}{\left(1-2\rho^{2}\cos \frac{2\pi j}{k} + \rho^{4}\right)^{\frac{N-2}{2}}}.
\end{equation}
Furthermore, since $A_{k}(\rho) \to -\infty$ as $\rho \to 0^{+}$ and $A_{k}(\rho) \to +\infty$ as $\rho \to 1^{-}$, there exists a unique $\rho_{0,k} \in (0,1)$ such that
\begin{equation}\label{imprho0}
A_{k}(\rho_{0,k})=0 \ \ \ \mbox{and} \ \ \ A_{k}(\rho)>0 \ \ \mbox{for all} \ \ \rho \in (\rho_{0,k},1).
\end{equation}
Now, in order to find the critical points of $\Psi_{k}$, we first fix $\rho\in(0,1)$ and require that the derivatives of $\Psi_{k}$ with respect to $d_0$ and $d_1$ vanish, namely
\begin{equation}\label{sistema}
\begin{cases}
d_0 + k \beta(\rho) d_1 - \displaystyle \frac{1}{d_0} =0, \\[2mm]
k A_{k}(\rho) d_1 + k \beta(\rho) d_0 -  \displaystyle \frac{k}{d_1}=0.
\end{cases}
\end{equation}
In the following lemma we solve the system \eqref{sistema}.
\begin{lemma}\label{lem:system}
The system \eqref{sistema} has a solution $(d_0,d_1)\in(0,+\infty)^2$ if and only if $\rho \in (\rho_{0,k}, 1)$. For each such $\rho$, the solution is unique and it is given by
\begin{equation}\label{solsis}
\bd d_1(\rho) = \frac{1}{\sqrt{A_{k}(\rho)+\beta(\rho)\mathfrak t_{k}(\rho)}} \qquad \mbox{and} \qquad \bd d_0(\rho) = \mathfrak t_{k}(\rho) \bd d_1(\rho),
\end{equation}
where
\begin{equation}\label{tkrhoimp}
\mathfrak t_{k}(\rho) : = \frac{\sqrt{(k-1)^{2}\beta^{2}(\rho)+4A_{k}(\rho)}-(k-1)\beta(\rho)}{2}.
\end{equation}
Moreover, the functions $\bd d_0(\rho)$ and $\bd d_1(\rho)$ are analytic on $(\rho_{0,k},1)$ and the matrix $D^{2}_{(d_0, d_1)} \Psi_{k}(\bd d_0(\rho), \bd d_1(\rho),\rho)$ is positive definite for every $\rho \in (\rho_{0,k}, 1)$.
\end{lemma}
\begin{proof}
First of all, we can rewrite the system \eqref{sistema} as
\begin{equation}\label{sistema1}
\begin{cases}
d_0^{2} + k \beta(\rho) d_0 d_1 -1=0, \\
A_{k}(\rho) d_1^{2} + \beta(\rho) d_0 d_1-1=0.
\end{cases}
\end{equation}
Subtracting the two equations in \eqref{sistema1} and dividing by $d_1^2$, we see that the ratio $t:=\frac{d_0}{d_1}>0$ satisfies
\begin{equation}\label{eqtk}
t^{2} + (k-1) \beta(\rho)\, t=A_{k}(\rho).
\end{equation}
Since $\beta(\rho)>0$, the equation \eqref{eqtk} has a positive solution if and only if $A_k(\rho)>0$, that is, by \eqref{imprho0} and the monotonicity of $A_k$, if and only if $\rho\in(\rho_{0,k},1)$. In this case the positive solution is unique and it is exactly $\mathfrak t_{k}(\rho)$ given in \eqref{tkrhoimp}; then the second equation in \eqref{sistema1} reads $$d_1^2\big(A_k(\rho)+\beta(\rho)\mathfrak t_k(\rho)\big)=1,$$ and \eqref{solsis} follows. Conversely, if $\rho\in(\rho_{0,k},1)$ and $\bd d_0(\rho)$, $\bd d_1(\rho)$ are given by \eqref{solsis}, then the second equation in \eqref{sistema1} holds by construction, and the first one follows from \eqref{eqtk}, since
$$\bd d_0^2+k\beta\,\bd d_0\bd d_1=\bd d_1^2\big(\mathfrak t_k^2+k\beta\,\mathfrak t_k\big)=\bd d_1^2\big(A_k+\beta\,\mathfrak t_k\big)=1.$$
Furthermore, by \eqref{betarho} and \eqref{Akk}, the functions $A_k$ and $\beta$ are analytic on $(0,1)$. Since $A_k(\rho)>0$ on $(\rho_{0,k},1)$ by \eqref{imprho0}, the radicand in \eqref{tkrhoimp} is positive and $\mathfrak t_k(\rho)>0$, so that $\mathfrak t_k$ is analytic on $(\rho_{0,k},1)$; moreover $A_k+\beta\,\mathfrak t_k>0$ there, and hence, by \eqref{solsis}, also $\bd d_1$ and $\bd d_0=\mathfrak t_k\bd d_1$ are analytic on $(\rho_{0,k},1)$.
Now, for every $\rho \in (\rho_{0,k},1)$, a direct computation shows that
$$D^{2}_{(d_0, d_1)} \Psi_{k}(\bd d_0(\rho), \bd d_1(\rho),\rho) =
\left(\begin{matrix}
1+\bd d_0^{-2}(\rho) & k \beta(\rho) \\
k \beta(\rho) & kA_{k}(\rho)+k\bd d_1^{-2}(\rho)
\end{matrix}\right)$$
and, using \eqref{sistema1} and \eqref{eqtk}, we get
$$D^{2}_{(d_0, d_1)} \Psi_{k}(\bd d_0(\rho), \bd d_1(\rho),\rho) =
\left(\begin{matrix}
2+\frac{k\beta(\rho)}{\mathfrak t_{k}(\rho)} & k \beta(\rho) \\
k \beta(\rho) & 2k A_{k}(\rho) + k\beta(\rho)\mathfrak t_{k}(\rho)
\end{matrix}\right).$$
Therefore, $\operatorname{tr} D^{2}_{(d_0, d_1)} \Psi_{k}(\bd d_0(\rho), \bd d_1(\rho),\rho)>0$ and
$$\det D^{2}_{(d_0, d_1)} \Psi_{k}(\bd d_0(\rho), \bd d_1(\rho),\rho) = 2k \left( 2A_{k}(\rho) + \beta(\rho) \mathfrak t_{k}(\rho) + \frac{k A_{k}(\rho) \beta(\rho)}{\mathfrak t_{k}(\rho)}\right)>0.$$
Thus this matrix is positive definite, and the claim follows.
\end{proof}
\noindent Now, let us define the function $f_{k} : (\rho_{0,k},1) \to \R$ by
\begin{equation}\label{fk}
\boxed{f_{k}(\rho) := \Psi_{k}(\bd d_0(\rho), \bd d_1(\rho), \rho) = \frac{k+1}{2} - \log \big[ \bd d_0(\rho) \bd d_1^{k}(\rho)\big],}
\end{equation}
where the second equality follows by summing the first equation in \eqref{sistema1} and $k$ times the second one.
The problem of finding the critical points of $\Psi_{k}$ can then be reduced to a one-dimensional problem, namely the study of the critical points of the function $f_{k}$.

The next result shows that the critical points of $\Psi_k$, together with their nondegeneracy and their Morse index, are completely determined by those of $f_k$; here the positive definiteness of $D^2_{(d_0,d_1)}\Psi_k$ at $(\bd d_0(\rho),\bd d_1(\rho),\rho)$ plays a crucial role.
\begin{proposition}\label{prop:link}
\begin{enumerate}[label=(\roman*)]
\item A point $(d_0,d_1,\rho)\in\mathcal M$ is a critical point of $\Psi_k$ if and only if $\rho\in(\rho_{0,k},1)$ is a critical point of $f_k$, $d_0=\bd d_0(\rho)$ and $d_1=\bd d_1(\rho)$.
\item Let $\bar\rho\in(\rho_{0,k},1)$ be a critical point of $f_k$ and let $\bar p:=(\bd d_0(\bar\rho),\bd d_1(\bar\rho),\bar\rho)$. Then
\begin{equation}\label{detSchur}
\det D^2\Psi_k(\bar p)=\det D^{2}_{(d_0, d_1)}\Psi_k(\bar p)\cdot f_k''(\bar\rho).
\end{equation}
In particular, $\bar p$ is a nondegenerate critical point of $\Psi_k$ if and only if $f_k''(\bar\rho)\neq0$. In this case, the Morse index of $\Psi_k$ at $\bar p$ is $0$ if $f_k''(\bar\rho)>0$ (and $\bar p$ is a strict local minimum point of $\Psi_k$), and it is $1$ if $f_k''(\bar\rho)<0$. Moreover, if $U\subset\mathcal M$ is an open neighbourhood of $\bar p$ such that $\overline U\subset\mathcal M$ and $\overline U $ contains no other critical point of $\Psi_k$, then
\begin{equation}\label{degree}
\deg\big(\nabla\Psi_k,U,0\big)=\sign f_k''(\bar\rho)\neq0.
\end{equation}
\end{enumerate}
\end{proposition}
\begin{proof}
$(i)$ If $(d_0,d_1,\rho)$ is a critical point of $\Psi_k$, then $(d_0,d_1)$ solves \eqref{sistema}, so by Lemma~\ref{lem:system} we have $\rho\in(\rho_{0,k},1)$, $d_0=\bd d_0(\rho)$ and $d_1=\bd d_1(\rho)$. On the other hand, for every $\rho\in(\rho_{0,k},1)$, since $\partial_{d_0}\Psi_k(\bd d_0(\rho),\bd d_1(\rho),\rho)=\partial_{d_1}\Psi_k(\bd d_0(\rho),\bd d_1(\rho),\rho)=0$, the chain rule gives
\begin{equation}\label{fkprimo-psi}
f_k'(\rho)=\partial_\rho\Psi_k\big(\bd d_0(\rho),\bd d_1(\rho),\rho\big).
\end{equation}
Hence $(\bd d_0(\rho),\bd d_1(\rho),\rho)$ is a critical point of $\Psi_k$ if and only if $f_k'(\rho)=0$.

$(ii)$ For brevity, we write $d=(d_0,d_1)$, $\bd d(\rho)=(\bd d_0(\rho),\bd d_1(\rho))$ and
\[
P:=D^2_{d}\Psi_k(\bar p),\qquad v:=\partial_\rho\nabla_{d}\Psi_k(\bar p)\in\R^2,\qquad s:=\partial^2_\rho\Psi_k(\bar p),
\]
so that
\[
D^2\Psi_k(\bar p)=\begin{pmatrix}P&v\\ v^{T}&s\end{pmatrix}.
\]
Differentiating the identity $\nabla_{d}\Psi_k(\bd d(\rho),\rho)=0$ at $\rho=\bar\rho$ we get $P\bd d'(\bar\rho)+v=0$, namely $$\bd d'(\bar\rho)=-P^{-1}v,$$ since $P$ is positive definite by Lemma~\ref{lem:system}. Then, differentiating \eqref{fkprimo-psi},
\[
f_k''(\bar\rho)=s+\langle v,\bd d'(\bar\rho)\rangle=s-\langle P^{-1}v,v\rangle,
\]
that is, $f_k''(\bar\rho)$ is the Schur complement of $P$ in $D^2\Psi_k(\bar p)$. Setting
\[
L:=\begin{pmatrix}I_2&P^{-1}v\\ 0&1\end{pmatrix},\qquad\text{we have}\qquad
D^2\Psi_k(\bar p)=L^{T}\begin{pmatrix}P&0\\ 0&f_k''(\bar\rho)\end{pmatrix}L,
\]
and \eqref{detSchur} follows since $\det L=1$. Moreover, by Sylvester's law of inertia, $D^2\Psi_k(\bar p)$ has the same number of negative (respectively zero) eigenvalues as $\mathrm{diag}\big(P,f_k''(\bar\rho)\big)$. Since $P$ is positive definite, the Morse index of $\Psi_k$ at $\bar p$ is $0$ if $f_k''(\bar\rho)>0$ and $1$ if $f_k''(\bar\rho)<0$; in the first case $D^2\Psi_k(\bar p)$ is positive definite, so that $\bar p$ is a strict local minimum point. Finally, if $f_k''(\bar\rho)\neq0$, then $\bar p$ is a nondegenerate zero of $\nabla\Psi_k$, hence
\[
\deg\big(\nabla\Psi_k,U,0\big)=\sign\det D^2\Psi_k(\bar p)=\sign f_k''(\bar\rho),
\]
where we used \eqref{detSchur} and $\det P>0$.
\end{proof}

\subsection{The derivative of \texorpdfstring{$f_{k}$}{f\_k}}
We now derive a convenient expression of the derivative of $f_{k}$. By the chain rule,
\begin{equation}
\begin{aligned}
f'_{k}(\rho) & = \frac{\partial \Psi_{k}(\bd d_0(\rho), \bd d_1(\rho), \rho)}{\partial d_0} \bd d_0'(\rho) + \frac{\partial \Psi_{k}(\bd d_0(\rho), \bd d_1(\rho), \rho)}{\partial d_1} \bd d_1'(\rho) + \frac{\partial \Psi_{k}(\bd d_0(\rho), \bd d_1(\rho), \rho)}{\partial \rho}
\\& = \frac{\partial \Psi_{k}(\bd d_0(\rho), \bd d_1(\rho), \rho)}{\partial \rho} = \frac{k}{2} \bd d_1^{2}(\rho) A'_{k}(\rho) + k \bd d_0(\rho) \bd d_1(\rho) \beta'(\rho),
\end{aligned}
\end{equation}
where in the second equality we used that $\frac{\partial \Psi_{k}}{\partial d_0}(\bd d_0(\rho), \bd d_1(\rho), \rho)=\frac{\partial \Psi_{k}}{\partial d_1}(\bd d_0(\rho), \bd d_1(\rho), \rho)=0$, since $(\bd d_0(\rho), \bd d_1(\rho))$ solves \eqref{sistema}.
Thus, since $\beta'(\rho)=-(N-2)\rho^{1-N}$ and $\bd d_0=\mathfrak t_k\bd d_1$, we obtain the following expression
\begin{equation}\label{fkprimoutile}
f'_{k}(\rho) = k(N-2) \rho^{1-N} \bd d_1^{2}(\rho) \left( \frac{\rho^{N-1}}{2(N-2)} A'_{k}(\rho) - \mathfrak t_{k}(\rho) \right).
\end{equation}
In Sections~\ref{sec:k3} and \ref{nE4} we will study the zeros of $f'_{k}$ in $(\rho_{0,k}, 1)$, and hence determine the critical points of the function $\Psi_{k}$.

\begin{lemma}\label{lem:limits}
We have
\begin{equation}\label{impfk}
\lim_{\rho\to\rho_{0,k}^{+}}f_k(\rho)=-\infty
\qquad\text{and}\qquad
\lim_{\rho\to1^{-}}f_k(\rho)=+\infty .
\end{equation}
\end{lemma}

\begin{proof}
As $\rho\to\rho_{0,k}^{+}$ we have $A_k(\rho)\to0^+$ and $\beta(\rho)\to\beta(\rho_{0,k})>0$, so that, by \eqref{tkrhoimp},
\[
\mathfrak t_k=\frac{2A_k}{\sqrt{(k-1)^2\beta^2+4A_k}+(k-1)\beta}=\frac{A_k}{(k-1)\beta}\,(1+o(1)),
\qquad
A_k+\beta\,\mathfrak t_k=\frac{k}{k-1}\,A_k\,(1+o(1)).
\]
Hence, by \eqref{solsis}, $\bd d_1=\big(\frac{k-1}{k A_k}\big)^{1/2}(1+o(1))$ and
$\bd d_0\bd d_1^{k}=\mathfrak t_k\bd d_1^{k+1}$ is of order $A_k^{\frac{1-k}{2}}\to+\infty$, so that $f_k\to-\infty$ by \eqref{fk}.

As $\rho\to1^{-}$ we have $A_k(\rho)\to+\infty$ and $\beta(\rho)\to0$, so that $\mathfrak t_k=\sqrt{A_k}\,(1+o(1))$ and $A_k+\beta\,\mathfrak t_k=A_k(1+o(1))$. Hence $\bd d_1=A_k^{-1/2}(1+o(1))\to0$ and $\bd d_0=\mathfrak t_k\bd d_1\to1$, so that $\bd d_0\bd d_1^k\to0$ and $f_k\to+\infty$ by \eqref{fk}.
\end{proof}

\section{The case \texorpdfstring{$k=3$}{k=3}}\label{sec:k3}
\noindent Throughout this section $k=3$ and we write $\rho_0:=\rho_{0,3}$. First of all we recall that, by \eqref{Akk},
\begin{equation}\label{A3}
A_{3}(\rho) = \frac1{(1-\rho^2)^{N-2}} -\frac2{(\sqrt3\rho)^{N-2}}+\frac2{(1+\rho^2+\rho^4)^{(N-2)/2}}.
\end{equation}
Then, we let
\begin{equation}\label{T}
T(\rho) : = \frac{\rho^{N-1}}{2(N-2)}A'_{3}(\rho)
=\frac{\rho^{N}}{(1-\rho^2)^{N-1}}+\frac1{3^{(N-2)/2}}-\frac{\rho^{N}(1+2\rho^2)}{(1+\rho^2+\rho^4)^{N/2}},
\end{equation}
where the second expression follows by differentiating \eqref{A3}, and
\begin{equation}\label{Q}
Q(\rho):=T^{2}(\rho) + 2 \beta(\rho) T(\rho)-A_{3}(\rho)=(T(\rho)-\mathfrak t_{3}(\rho))(T(\rho)+\mathfrak t_{3}(\rho)+2\beta(\rho)).
\end{equation}
Now, combining \eqref{fkprimoutile} with the fact that $A'_{3}(\rho), \mathfrak t_{3}(\rho)$ and $\beta(\rho)$ are all positive on $(\rho_0,1)$, we observe that
\begin{equation}\label{signo}
\sign f'_{3}(\rho) = \sign \big(T(\rho)-\mathfrak t_{3}(\rho)\big) = \sign Q(\rho) \ \ \ \ \mbox{on} \ \ (\rho_{0},1).
\end{equation}
More precisely, by \eqref{fkprimoutile} and \eqref{Q},
\begin{equation}\label{signog}
f_3'(\rho)=g(\rho)\,Q(\rho),\qquad g(\rho):=\frac{3(N-2)\rho^{1-N}\bd d_1^2(\rho)}{T(\rho)+\mathfrak t_3(\rho)+2\beta(\rho)}>0,\qquad \rho\in(\rho_0,1),
\end{equation}
so that, since $g$ is smooth, $f_3''(\rho)=g(\rho)Q'(\rho)$ at every zero $\rho$ of $Q$.

The next lemma will be useful in the following.
\begin{lemma}\label{Pos}
Let $N\geq4$. Then $Q(\rho)>0$ for every $\rho\in\left[\frac{1}{\sqrt3},1\right)$. In particular, $f_3'(\rho)>0$ for every $\rho\in(\rho_0,1)$ such that $\rho\geq1/\sqrt3$.
\end{lemma}

\begin{proof}
First of all, using the definitions of the functions $T$ and $Q$ as in \eqref{T} and \eqref{Q}, we can define
\begin{equation}\label{FN}
\mathcal F_N(\rho)
: =(\sqrt3\rho)^{N-2}\bigl(Q(\rho)-T^{2}(\rho)\bigr)=(\sqrt3\rho)^{N-2}\left[
\frac{\rho(1-\rho^{N-2})}{N-2}A_3'(\rho)-A_3(\rho)
\right].
\end{equation}
Therefore
\begin{equation}\label{utileQrho}
Q(\rho)=T(\rho)^2+\frac{\mathcal F_N(\rho)}{(\sqrt3\rho)^{N-2}}
\end{equation}
and so, in order to get the claim, it suffices to show that $\mathcal F_N(\rho)>0$ for $\frac{1}{\sqrt3}\leq\rho<1$. Now, we set
\begin{equation}\label{ABC}
\mathfrak a=\frac{3\rho^2}{(1-\rho^2)^2},\qquad
\mathfrak b=\frac{3\rho^2}{1+\rho^2+\rho^4},\qquad
\mathfrak c=\frac{\rho^2(1+2\rho^2)}{1+\rho^2+\rho^4}.
\end{equation}
A straightforward computation shows that $0<\mathfrak b <1$ and $0<\mathfrak c<1$ for all $0<\rho<1$. Moreover, by \eqref{A3} and \eqref{T}, letting
\begin{equation}\label{D}
\mathfrak d: =\frac{2\rho^2}{1-\rho^2},
\end{equation}
we can rewrite \eqref{FN} as
\begin{equation}\label{eqtriF}
\mathcal F_N(\rho)
=\mathfrak a^{\frac{N-2}{2}}\left[\mathfrak d(1-\rho^{N-2})-1\right]
+2\left(1-\mathfrak b^{\frac{N-2}{2}}\right)+2(1-\rho^{N-2})\left(1-\mathfrak c \mathfrak b^{\frac{N-2}{2}}\right).
\end{equation}
Note that $\mathfrak d\geq1$ if and only if $\rho\geq\frac1{\sqrt3}$.
To prove the positivity of $\mathcal F_N$, we distinguish two ranges of $\rho$.

\smallskip
\noindent
\emph{Case 1: $1/\sqrt2\leq\rho<1$.} Since $\rho^{N-2}\leq\rho^2$, we obtain that $\mathfrak d(1-\rho^{N-2})-1\geq\mathfrak d(1-\rho^{2})-1=2\rho^2-1\geq0$. Moreover, since $0<\mathfrak b<1$ and $0<\mathfrak c<1$, the last two terms of \eqref{eqtriF} are strictly positive; hence $\mathcal F_N(\rho)>0$ in this range.\\

\smallskip
\noindent
\emph{Case 2: $1/\sqrt3\leq\rho\leq1/\sqrt2$.} Recalling the definition of $\mathfrak a$ and $\mathfrak d$ as in \eqref{ABC} and \eqref{D}, we can observe that
\begin{equation*}
(\mathfrak d-1)\mathfrak a+\frac34-\mathfrak d\rho^2\mathfrak a=\mathfrak a\bigl[\mathfrak d(1-\rho^2)-1\bigr]+\frac34=\mathfrak a(2\rho^2-1)+\frac34=\frac{3(1-3\rho^2)^2}{4(1-\rho^2)^2}.
\end{equation*}
This means that
$$(\mathfrak d-1)\mathfrak a+\frac34-\mathfrak d\rho^2\mathfrak a \geq 0$$
and so
\begin{equation*}
\rho^2\mathfrak a\leq\frac{\mathfrak d-1}{\mathfrak d}\mathfrak a+\frac{1}{\mathfrak d}\cdot\frac34.
\end{equation*}
Therefore, thanks to the convexity of the function $z\mapsto z^{\frac{N-2}{2}}$, we get that
\begin{equation*}
(\rho^2\mathfrak a)^{\frac{N-2}{2}}
\leq\left(\frac{\mathfrak d-1}{\mathfrak d}\mathfrak a+\frac{1}{\mathfrak d}\cdot\frac34\right)^{\frac{N-2}{2}}\leq\frac{\mathfrak d-1}{\mathfrak d}\mathfrak a^{\frac{N-2}{2}}+\frac{1}{\mathfrak d}\left(\frac34\right)^{\frac{N-2}{2}}.
\end{equation*}
It follows that
\begin{equation*}
\mathfrak a^{\frac{N-2}{2}}\bigl[\mathfrak d(1-\rho^{N-2})-1\bigr]
=(\mathfrak d-1)\mathfrak a^{\frac{N-2}{2}}-\mathfrak d(\rho^2\mathfrak a)^{\frac{N-2}{2}}\geq-\left(\frac34\right)^{\frac{N-2}{2}}\geq-\frac34.
\end{equation*}
Moreover, with a straightforward computation one can show that for all $\frac{1}{\sqrt3} \leq \rho \leq \frac{1}{\sqrt2}$,
$$\mathfrak b\leq\frac{6}{7} \ \ \ \mbox{and} \ \ \ \mathfrak c \leq \frac{4}{7}.$$
Therefore, using the fact that $\mathfrak b^{\frac{N-2}{2}} \leq \mathfrak b$ and $\rho^{N-2}\leq\rho^2\leq1/2$, from \eqref{eqtriF} we obtain that
\begin{equation*}
\mathcal F_N(\rho)
\geq-\frac34+2(1-\mathfrak b)+(1-\mathfrak c\mathfrak b)\geq-\frac34+\frac27+\frac{25}{49}>0.
\end{equation*}
Finally, combining \emph{Case 1} and \emph{Case 2}, we get that $\mathcal F_N(\rho)>0$ and hence $Q(\rho)>0$, for every $\rho\in\left[\frac{1}{\sqrt3},1\right)$.
\end{proof}
\begin{remark}\label{rmk1}
By Lemma~\ref{Pos} and \eqref{signo}, for every
$N\geq4$, any critical point $\rho\in(\rho_0,1)$ of $f_3$
satisfies the condition $\rho<\frac{1}{\sqrt3}$.
\end{remark}

\subsection{Existence for \texorpdfstring{$N \geq 13$}{N>=13}}
We first show the existence of exactly two critical points of $f_3$.
\begin{proposition}\label{3.1}
Let $N\geq13$. The function $f_3$ has exactly two critical points in $(\rho_0,1)$,
both nondegenerate: a local maximum $\rho^*$ and a local minimum $\rho^{**}$. Moreover,
\begin{equation}\label{critrad}
\rho_0<\rho^*<\frac{11}{20}<\rho^{**}<\frac{1}{\sqrt3},
\end{equation}
and
\[
f_3''(\rho^*)<0<f_3''(\rho^{**}).
\]
\end{proposition}

\begin{proof}
\emph{Step 1: existence of two critical points.} In order to prove that the function $f_{3}$ admits at least two critical points it \del{is} suffices to show that $r_{*}:=\frac{11}{20}$ satisfies the two following conditions
\begin{equation}\label{condimp!}
A_{3}(r_{*})>0 \ \ \ \ \mbox{and} \ \ \ f_{3}'(r_{*})<0.
\end{equation}
Indeed, for all $\rho \in (0,1)$, we know that $A'_{3}(\rho)>0$ and so, $A_{3}(\rho)$ is strictly increasing. Therefore, since $A_{3}(\rho_{0})=0$, the first inequality of \eqref{condimp!} ensures that $r_{*} \in (\rho_{0},1)$. \ \\
Furthermore, the second inequality of \eqref{condimp!} combined with \eqref{impfk} implies that the function $f_{3}$ has a local maximum at some $\rho^{*} \in (\rho_{0}, r_{*})$ and a local minimum at some $\rho^{**} \in (r_{*},1)$. \ \\
Therefore, in order to get the claim we have to prove that the two inequalities in \eqref{condimp!} hold. \ \\
Let us set the following quantity
\begin{equation}
\mathfrak R := (1-r_{*}^{2})^{2-N}\end{equation}
Using \eqref{A3}, we easily obtain that
\begin{equation}\label{A3rstar}
A_{3}(r_*)>\mathfrak R\left[1-2\left(\frac{1-r_*^2}{\sqrt3\,r_*}\right)^{N-2}\right]>0
\end{equation}
and so the first inequality in \eqref{condimp!} is verified. \ \\
Now, in order to prove the second inequality in \eqref{condimp!}, using \eqref{signo}, we can just verify that $Q(r_{*})<0$. \ \\
First of all, we differentiate \eqref{A3} with respect to $\rho$, obtaining $T(\rho)$ defined in \eqref{T}.
Hence, with a direct computation, dropping the negative last term in \eqref{T}, one can show that
\begin{equation*}
0<\frac{T(r_*)}{r_*^{N-2}\mathfrak R}<\frac{r_*^2}{1-r_*^2}+\left(\frac{1-r_*^2}{\sqrt3\,r_*}\right)^{N-2}<1.
\end{equation*}
Now, noting that
$$r_*^{2(N-2)}\mathfrak R<\frac{1}{2000} \ \ \mbox{and} \ \ \beta(r_*)<r_*^{2-N},$$
we obtain that
\begin{equation}\label{stime}
\frac{T^{2}(r_*)}{\mathfrak R}<\frac1{2000},\qquad \frac{2\beta(r_*)T(r_*)}{\mathfrak R}<\frac{242}{279}+\frac{33}{500}.
\end{equation}
Finally, combining the estimates \eqref{stime} with the fact that $\frac{A_{3}(r_{*})}{\mathfrak R} > \frac{467}{500}$, using the definition of $Q(\rho)$ given in \eqref{Q}, we easily get that $\frac{Q(r_*)}{\mathfrak R}<0$ and so the claim follows. \ \\ \\
\emph{Step 2: uniqueness of the two critical points.}
We want to show that the function $Q$, and hence $f_3'$, has at most two zeros in $(\rho_0,1)$.
Since $Q>0$ on $\big[\frac1{\sqrt3},1\big)$ by Lemma~\ref{Pos}, it is enough to show that $Q$ has at most two zeros in $\left(0, \frac{1}{\sqrt3}\right)$ and this follows once we know that $Q'$ has at most one zero there. Indeed, suppose that $Q$ vanishes at
three different points $z_1<z_2<z_3$ of $\left(0,\frac{1}{\sqrt3}\right)$. Since $Q$ is a smooth function, Rolle's theorem on $[z_1,z_2]$ and on $[z_2,z_3]$ gives
$\bar{z}_1\in(z_1,z_2)$ and $\bar{z}_2\in(z_2,z_3)$ with $Q'(\bar{z}_1)=Q'(\bar{z}_2)=0$. These are two distinct zeros of $Q'$ in $\left(0,\frac{1}{\sqrt3}\right)$, which is a contradiction. \ \\
Hence, it remains to prove that $Q'$ has at most one zero in
$\left(0,\frac{1}{\sqrt{3}}\right)$.
To this end, we use the expansion
\begin{equation}\label{powerseries}
(1-\rho^2)^{-(N-2)}
+2(1+\rho^2+\rho^4)^{-\frac{N-2}{2}}
=\sum_{j=0}^{\infty}c_j\rho^{2j},
\qquad |\rho|<1,
\end{equation}
where $c_j\geq0$ for every $j\geq0$ and
\[
c_0=3,\qquad c_1=0,\qquad
c_2=\frac{3(N-2)^2}{4}.
\]
These properties are established in Lemma~\ref{lem:positive-coefficients} in the Appendix. \ \\
Now, by \eqref{A3} and \eqref{powerseries}, we get
\begin{equation}\label{eq:tri-AT-series}
A_{3}(\rho)=\sum_{j\geq0}c_j\rho^{2j}-2\,3^{-\frac{N-2}{2}}\rho^{2-N},\qquad
T(\rho)=\frac{\rho^{N-1}}{2(N-2)}A'_{3}(\rho)=3^{-\frac{N-2}{2}}+S(\rho),
\end{equation}
where
\[
S(\rho):=\frac1{N-2}\sum_{j\geq2}jc_j\rho^{N-2+2j}.
\] 
Since $c_j\geq0$ and $c_2>0$, we have that $S(\rho)>0$ and $T'(\rho)=S'(\rho)>0$ for all $\rho \in (0,1)$.
\ \\
Furthermore, using
\[
A_3'(\rho)=2(N-2)\rho^{1-N}T(\rho),
\qquad
\beta'(\rho)=-(N-2)\rho^{1-N}
\quad \mbox{and} \ \ \rho^{N-2}\beta(\rho)=1-\rho^{N-2}, \] 
and computing the first derivative of the function \eqref{Q} we obtain that
\begin{equation}\label{Qprimo}
Q'(\rho)=2T'(\rho)\bigl(T(\rho)+\beta(\rho)\bigr)
-4(N-2)\rho^{1-N}T(\rho)=\underbrace{2\rho^{2-N}T'(\rho)}_{(>0)}\Theta(\rho),
\end{equation}
where
\[
\Theta(\rho):=
1-\rho^{N-2}\bigl(1-T(\rho)\bigr)
-2(N-2)\frac{T(\rho)}{\rho T'(\rho)}.
\]
Therefore, applying Lemma~\ref{lem:auxiliary-monotonicity} with $p=N-2$ and $a=3^{-(N-2)/2}$ (its assumption $\rho^{p}<8a$ reads $(\sqrt3\rho)^{N-2}<8$, which holds for every $\rho<\frac1{\sqrt3}$), it follows that
$\Theta'(\rho)>0$ for $0<\rho<\frac{1}{\sqrt{3}}$. \ \\
Moreover, since \eqref{Qprimo} holds, arguing by contradiction one can easily show that $Q'$ has at most one zero in this interval. \ \\
Finally, combining \emph{Step 1} and \emph{Step 2}, we can say that the function $Q$ has exactly two zeros in $(\rho_0,1)$, namely $\rho^*$ and $\rho^{**}$, and by Lemma \ref{Pos} we have $\rho^{**}<\frac{1}{\sqrt3}$, so that \eqref{critrad} holds. Since $Q(\rho^*)=Q(\rho^{**})=0$, by Rolle's theorem $Q'$ has a zero in $(\rho^*,\rho^{**})\subset\left(0,\frac1{\sqrt3}\right)$. Since $\Theta$ is strictly increasing on $\left(0,\frac1{\sqrt3}\right)$, by \eqref{Qprimo} we get that $Q'$ is negative on the left of this zero and positive on the right. Therefore $Q'(\rho^*)<0<Q'(\rho^{**})$, and hence, by \eqref{signog}, $f_3''(\rho^*)<0<f_3''(\rho^{**})$.
\ \\ \\

\end{proof}

\subsection{Nonexistence for \texorpdfstring{$3\leq N\leq12$}{3<=N<=12}}
In the following lemma we show that $f_3$ has no critical points for $3\leq N\leq12$. The proof is different for $4\leq N\leq12$ and for $N=3$.

\begin{lemma}\label{nonexgeq4}
If $3\leq N\leq12$, then $f_3'(\rho)>0$ for every $\rho\in(\rho_0,1)$.
\end{lemma}

\begin{proof}
\emph{Case 1: $4\leq N\leq12$.} By Lemma~\ref{Pos}, we already know that $f_3'(\rho)>0$ whenever $\rho\geq\frac{1}{\sqrt3}$ and so it remains to consider the range $0<\rho<\frac{1}{\sqrt3}$.

Let $\mathcal F_N$, $\mathfrak a$, $\mathfrak b$ and $\mathfrak c$ be as in the proof of Lemma~\ref{Pos} (see \eqref{FN} and \eqref{ABC}).

By \eqref{utileQrho} and \eqref{signo}, in order to get the claim it suffices to show that
$$\mathcal F_N(\rho)>0 \ \ \mbox{for all} \ \ 0<\rho<\frac1{\sqrt3} \ \ \mbox{and} \ \ 4\leq N\leq12.$$
To this end, we let
$$\mathcal G(N,\rho):=\mathcal F_N(\rho)-2\mathfrak c(\rho^2\mathfrak b)^{\frac{N-2}{2}}.$$
Therefore, from \eqref{eqtriF}, we obtain that
\begin{equation}\label{eqtriFexp}
\mathcal G(N,\rho)
=4-2\rho^{N-2}-2(1+\mathfrak c)\mathfrak b^{\frac{N-2}{2}}-\frac{1-3\rho^2}{1-\rho^2}\mathfrak a^{\frac{N-2}{2}}
-\frac{2\rho^2}{1-\rho^2}(\rho^2\mathfrak a)^{\frac{N-2}{2}},
\end{equation}
where we regard $N$ as a real parameter in $[4,12]$.

Let us fix $\rho\in\left(0,\frac1{\sqrt3}\right]$. Since $1-3\rho^2\geq0$ and $\mathfrak c>0$, all the coefficients of the powers in \eqref{eqtriFexp} are nonpositive; hence the convexity of $N\mapsto u^{\frac{N-2}{2}}$ for $u>0$ implies that $\mathcal G(\cdot,\rho)$ is concave on $[4,12]$. Therefore, since $\mathfrak c>0$,
\begin{equation}\label{eqtriint}
\mathcal F_N(\rho)\geq\mathcal G(N,\rho)
\geq\min\{\mathcal G(4,\rho),\mathcal G(12,\rho)\}\qquad\text{for every }N\in[4,12].
\end{equation}
It remains to prove that $\mathcal G(4,\rho)>0$ and $\mathcal G(12,\rho)>0$.

In particular,
\begin{equation}
\begin{aligned}
\mathcal G(4,\rho) & = 4-2\rho^{2}-2(1+\mathfrak c) \mathfrak b - \left[ \frac{1-3\rho^{2}}{1-\rho^{2}} \mathfrak a + \frac{2\rho^{2}}{1-\rho^{2}} (\rho^{2}\mathfrak a) \right]
\\& =4-2\rho^2-2(1+\mathfrak c)\mathfrak b-\mathfrak a(1-2\rho^2).
\end{aligned}
\end{equation}
Moreover, for $0<\rho\leq\frac1{\sqrt3}$, we have
\begin{equation}\label{boundbc}
\mathfrak b\leq\frac9{13} \ \ \ \mbox{and} \ \ \ \mathfrak c\leq\frac{5}{13},
\end{equation}
so that $2(1+\mathfrak c)\mathfrak b<2$.

Furthermore, a straightforward computation shows that $\mathfrak a(1-2\rho^{2}) \leq \frac 34$, and so
\begin{equation}\label{eqtriG4}
\mathcal G(4,\rho)
=4-2\rho^2-2(1+\mathfrak c)\mathfrak b-\mathfrak a(1-2\rho^2)>4-\frac23-2-\frac34=\frac7{12}>0.
\end{equation}
Now, we will prove that $\mathcal G(12,\rho)>0$. By \eqref{eqtriFexp}, we have
\begin{equation*}
\mathcal G(12,\rho)=4-2\rho^{10}-2(1+\mathfrak c)\mathfrak b^5-\frac{1-3\rho^2}{1-\rho^2}\mathfrak a^5-\frac{2\rho^2}{1-\rho^2}(\rho^2\mathfrak a)^5.
\end{equation*}
Therefore, combining the bounds on $\mathfrak b$ and $\mathfrak c$ given in \eqref{boundbc}, with
$$\rho^2\leq\frac13, \ \ \ \ \ \ \frac{2\rho^2}{1-\rho^2}\leq1 \ \ \ \ \mbox{and} \ \ \ \rho^2\mathfrak a\leq\frac34,$$
we obtain
\begin{equation}\label{lowG12a}
2\rho^{10}+2(1+\mathfrak c)\mathfrak b^5+\frac{2\rho^2}{1-\rho^2}(\rho^2\mathfrak a)^5\leq\frac2{243}+\frac{36}{13}\left(\frac9{13}\right)^5+\left(\frac34\right)^5<\frac34.
\end{equation}
Moreover, since $1-3\rho^2\geq0$, the weighted arithmetic-geometric mean inequality gives
\begin{equation*}
(3\rho^2)^{\frac5{11}}\cdot \bigl[10(1-3\rho^2)\bigr]^{\frac1{11}}
\leq\frac5{11}(3\rho^2)+\frac5{11}
+\frac1{11}\cdot10(1-3\rho^2)
=\frac{15}{11}(1-\rho^2),
\end{equation*}
and so we get
\begin{equation}\label{lowG12b}
\frac{1-3\rho^2}{1-\rho^2}\mathfrak a^5
=\frac{3^5\rho^{10}(1-3\rho^2)}{(1-\rho^2)^{11}}
\leq\frac1{10}\left(\frac{15}{11}\right)^{11}<\frac{13}{4}.
\end{equation}
By \eqref{lowG12a} and \eqref{lowG12b}, we conclude that $\mathcal G(12,\rho)>4-\frac34-\frac{13}{4}=0$.

Finally, combining \eqref{eqtriG4} with the positivity of $\mathcal G(12,\rho)$ and \eqref{eqtriint}, we obtain that $\mathcal F_N(\rho)>0$ for every $0<\rho\leq\frac1{\sqrt3}$ and $4\leq N\leq12$.

Hence $Q>0$, and so $f_3'>0$, in $\left(\rho_0,\frac1{\sqrt3}\right]$. Together with Lemma~\ref{Pos}, we obtain $f_3'>0$ on $(\rho_0,1)$, which proves the claim for $4\leq N\leq12$.

\smallskip
\noindent
\emph{Case 2: $N=3$.} By \eqref{signo}, in order to prove the claim it suffices to show that $Q>0$ on $(0,1)$.

Now, we recall the expression of the coefficients of the expansion as in Lemma~\ref{lem:positive-coefficients}, namely
\[
c_0=3,\qquad c_1=0,\qquad c_2=\frac34,\qquad
c_3=\frac{15}{8},
\]
with $c_j\geq0$ for every $j$.

Moreover, using \eqref{eq:tri-AT-series} and $T^2\geq \frac{2T}{\sqrt3}-\frac13$ in $Q=T^2+2\beta T-A$, we obtain
\[
Q(\rho)\geq
\frac1{\sqrt3}\left(\frac4\rho-2\right)-\frac83
+\sum_{j\geq2}c_j\rho^{2j}
\left[2j-1-2j\left(1-\frac1{\sqrt3}\right)\rho\right].
\]
First we note that, for all $0<\rho<1$, since $\frac{1}{\sqrt3}>\frac{4}{7}$,
\begin{equation}
2j-1-2j\left(1-\frac1{\sqrt3}\right)\rho > \frac{8}{7}j-1>0.
\end{equation}
Keeping only the terms with $j=2,3$, we deduce
\[
\begin{aligned}
Q(\rho)
&>\frac{16}{7\rho}-\frac{80}{21}
+\frac{27}{28}\rho^4+\frac{255}{56}\rho^6>\frac{16}{21\rho}\left(3-5\rho+7\rho^7\right).
\end{aligned}
\]
Finally, the weighted arithmetic--geometric mean inequality yields
\[
3+7\rho^7
=7\left(\frac67\cdot\frac12+\frac17\cdot7\rho^7\right)
\geq7\left(\frac7{64}\right)^{1/7}\rho
>5\rho.
\]
Thus $Q>0$ on $(0,1)$, and the claim for $N=3$ follows from
\eqref{signo}.
\end{proof}

\subsection{Proof of Theorem~\ref{thm:main}}

\begin{proof}[Proof of Theorem~\ref{thm:main}]
$(i)$ Let $N\geq13$. By Proposition~\ref{3.1} and
Proposition~\ref{prop:link}, $\Psi_3$ has exactly two
critical points in $\mathcal M$,
\[
p^*=\bigl(\bd d_0(\rho^*),\bd d_1(\rho^*),\rho^*\bigr),
\qquad
p^{**}=\bigl(\bd d_0(\rho^{**}),\bd d_1(\rho^{**}),
\rho^{**}\bigr),
\]
where $\rho^*<\rho^{**}$ satisfy \eqref{critrad}.
Both points are nondegenerate: $p^*$ has Morse index $1$,
whereas $p^{**}$ is a strict local minimum.
\\We now apply the standard degree argument for the reduced
functional, as in \cite{PRV3}.
Indeed, by \eqref{degree}, the local degrees of
$\nabla\Psi_3$ at $p^*$ and $p^{**}$ are $-1$ and $1$,
respectively, while \eqref{expansion} gives
\[
(\omega_N\e)^{-1}\nabla\widetilde J_\e
=\nabla\Psi_3+o(1)
\qquad\text{locally uniformly in }\mathcal M.
\]
Hence, for $\e$ small, $\widetilde J_\e$ has critical points
$p_\e^{(1)}$ and $p_\e^{(2)}$ in fixed disjoint compact
neighbourhoods of $p^*$ and $p^{**}$.
Every accumulation point of these critical points is a zero
of $\nabla\Psi_3$, so
\[
p_\e^{(1)}\to p^*,
\qquad
p_\e^{(2)}\to p^{**}.
\]
Therefore, Proposition~\ref{prop:reduction} yields the corresponding
solutions
\[
u_\e^{(i)}
=\mathcal V_\e(p_\e^{(i)})+\phi_\e(p_\e^{(i)}),
\qquad i=1,2,
\]
with $\|\phi_\e(p_\e^{(i)})\|\to0$.
Their concentration profiles show that they change sign
and are distinct, since their limiting concentration radii
are different.

$(ii)$ Let $3\leq N\leq12$.
By Lemma~\ref{nonexgeq4}, $f_3'>0$
on $(\rho_0,1)$; hence Proposition~\ref{prop:link}$(i)$
implies that $\Psi_3$ has no critical points in $\mathcal M$.
Therefore, by Lemma~\ref{lem:limit-config}, there is no family of solutions as in the statement.
\end{proof}

\begin{remark}[Limit configurations]\label{rmk:limit}
Let $N\geq13$ and let $p^*,p^{**}$ be as in the proof of Theorem~\ref{thm:main}. Then $p^*$ and $p^{**}$ are the only possible limit configurations of solutions of the form \eqref{ansatz}: if $\e_n\to0$ and $u_n=\mathcal V_{\e_n}+\phi_n\in\mathcal H_3$ are solutions of \eqref{pb} with $\e=\e_n$, $\phi_n\in\mathcal K^\perp_{d_{0,n},d_{1,n},\rho_n}$, $\|\phi_n\|\to0$ and $(d_{0,n},d_{1,n},\rho_n)\to p\in\mathcal M$, then $p\in\{p^*,p^{**}\}$. Indeed, by Lemma~\ref{lem:limit-config}, $p$ is a critical point of $\Psi_3$, and $p^*,p^{**}$ are the only critical points of $\Psi_3$ in $\mathcal M$.
\end{remark}

\begin{remark}[The case $k=2$: on a conjecture of \cite{BDP}]\label{rmk:k2}
The method of this section also settles the case of two antipodal negative bubbles, which was left open in \cite{BDP}. Lemma~\ref{lem:system}, Proposition~\ref{prop:link} and \eqref{fkprimoutile} do not use the assumption $k\geq3$, and hold for $k=2$, with $\xi_{1,2}(\rho)=(\pm\rho,0,\dots,0)$ (and with $\Gamma_2$, $\mathcal H_2$ defined as in \eqref{Hk}). Moreover $A_2$ coincides with the function $\alpha$ in \cite[Proposition~2.3]{BDP}, and the change of variables $\lambda=\sqrt{c_N/2}\,d_0$, $\mu=\sqrt{c_N/2}\,d_1$ gives $F=c_N\Psi_2+\mathrm{const}$, where $F$ is the reduced function in \cite[(2.3)]{BDP}. In \cite[Lemma~3.1]{BDP} it is proved that $f_2$ has a local maximum point $\rho_1<\frac12$ and a local minimum point $\rho_2>\frac12$, and in \cite[Remark~3.2]{BDP} it is conjectured that these are its only critical points. We claim that, \\\begin{center} \emph{for every $N\geq3$, $f_2$ has exactly two critical points $\rho_1<\frac12<\rho_2$ in $(\rho_{0,2},1)$, and $f_2''(\rho_1)<0<f_2''(\rho_2)$}. \end{center}By Proposition~\ref{prop:link}, $\Psi_2$ (equivalently, $F$) has exactly two critical points, a nondegenerate saddle point of Morse index $1$ and a nondegenerate strict local minimum, which proves the conjecture; in particular, by \cite[Theorem~2.6(b)]{BDP}, the two families of solutions found in \cite{BDP} are $C^1$ curves in every dimension.

Let $a:=2^{1-N}$ and $\bar\rho:=\frac{\sqrt5-1}{2}$. By \eqref{Akk},
\[
A_2(\rho)=\frac1{(1-\rho^2)^{N-2}}-\frac1{(2\rho)^{N-2}}+\frac1{(1+\rho^2)^{N-2}},\quad
T(\rho):=\frac{\rho^{N-1}}{2(N-2)}A_2'(\rho)=\frac{\rho^N}{(1-\rho^2)^{N-1}}+a-\frac{\rho^N}{(1+\rho^2)^{N-1}},
\]
and, let, as in \eqref{Q}, $Q=(T-\mathfrak t_2)(T+\mathfrak t_2+\beta)$. Exactly as in \eqref{signog} we have $f_2'=g\,Q$ with $g>0$ on $(\rho_{0,2},1)$. Moreover, $Q(\rho_{0,2})>0$, $Q(1/2)<0$ by \cite[(3.5)]{BDP},
and $Q(\rho)\to+\infty$ as $\rho\to1^-$. Since $T,\beta>0$, the inequality $Q(1/2)<0$ gives $A_2(1/2)>T^2(1/2)+\beta(1/2)T(1/2)>0$, so that $\rho_{0,2}<\frac12$.\\\\ Hence it suffices to show that $Q$ has exactly two zeros $\rho_1<\frac12<\rho_2$ in $(\rho_{0,2},1)$, with $Q'(\rho_1)<0<Q'(\rho_2)$.
\\The binomial expansion gives $T=a+S$, where $S$ is a power
series with nonnegative coefficients and $S'>0$. Differentiating,
we obtain
\[
Q'=\rho^{2-N}T'\Theta_2,\qquad
\Theta_2=1-\rho^{N-2}(1-2T)-3(N-2)\frac{T}{\rho T'}.
\]
By the Cauchy--Schwarz estimate in the proof of Lemma~\ref{lem:auxiliary-monotonicity},
\[
\rho\Theta_2'
\ge \frac{N-2}{S}\bigl(3a-\rho^{N-2} S(1-2T)\bigr)
     +2\rho^{N-1}T'.
\]
Since $S(1-2T)\le S(1-2S)\le1/8$, it follows that
$\Theta_2'>0$ wherever $(2\rho)^{N-2}<12$ or $T\ge1/2$.
Also, $T$ is increasing and
\[
\tfrac12<T(\bar\rho)
=a+\bar\rho\bigl(1-5^{-(N-1)/2}\bigr)<1.
\]
Thus $\Theta_2$ is strictly increasing on $(0,1/2]$ and on
$[\bar\rho,1)$. \\ Now, if $N\le7$, then $\Theta_2$ is strictly increasing on all
of $(0,1)$, since $(2\bar\rho)^5<12$.
\\ Instead, if $N\geq 8$, $Q$ has no zeros in  $[1/2,\bar\rho]$ since on this interval $0<T<1$ and, by using the fact that $\beta=\rho^{2-N}-1$ and the
explicit formulas for $T$ and $A_2$, we obtain
\[
Q<\rho^{2-N}T-A_2
\le \frac32(2\rho)^{2-N}
   -\frac{1-2\rho^2}{(1-\rho^2)^{N-1}}
\le \frac32-(1-2\bar\rho^2)\Bigl(\frac43\Bigr)^7<0.
\]

Since $Q'$ has the sign of $\Theta_2$, on each interval where
$\Theta_2$ is strictly increasing, $Q'$ can change sign only
once, from negative to positive. The above signs of $Q$ therefore
give exactly two zeros $\rho_1<1/2<\rho_2$, with
$Q'(\rho_1)<0<Q'(\rho_2)$. Hence
$f_2''(\rho_1)<0<f_2''(\rho_2)$, as claimed.

\smallskip
\end{remark}
 \section{The nodal sets}\label{nodal}
\noindent First, we describe the asymptotic behaviour of the solutions away from the concentration points. In particular, let us consider, for $b\in(0,1)$,
$$\mathcal A_b:=\{x\in\mathcal B:b<|x|<1\}.$$

\begin{lemma}\label{lem:nodal-farfield}
Let $u_\e$ be either family in Theorem~\ref{thm:main}, with
limiting parameters
$(\bd d_0(\rho),\bd d_1(\rho),\rho)$.
For every $b\in(\rho,1)$,
\begin{equation}\label{eq:nodal-farfield}
\frac{u_\e}{\sqrt\e}\longrightarrow\Phi_\rho
\quad\text{in }C^1(\overline{\mathcal A_b}),
\end{equation}
where
\begin{equation}\label{eq:nodal-profile}
\Phi_\rho(x)=\alpha_N\sqrt{\mathtt d_N}
\left[\bd d_0(\rho)G(x,0)
-\bd d_1(\rho)\sum_{j=1}^3G(x,\xi_j(\rho))\right].
\end{equation}
\end{lemma}

\begin{proof}
In order to establish the behaviour of the solutions away from the concentration points we argue as in \cite[Theorem 2.6(a)(iv)]{BDP} and \cite[Lemmas~4.3--4.5]{BMP}, applying their results to the present four-bubble configuration.
\\In particular, the concentration points stay mutually separated and away from the boundary and the parameters converge in $\mathcal M$. Moreover, the estimates away from the concentration points, together with the Green representation and standard boundary elliptic estimates (see for instance \cite[Theorem~9.13]{GT}), give the $C^1$ convergence on $\overline{\mathcal A_b}$ to the signed sum of Green functions.
\\Furthermore, for the normalization used to define the Green function in \eqref{Green}, the contribution of a projected bubble is equal to $\alpha_N\delta^{(N-2)/2}G(x,\xi)$; since \eqref{deltas} gives
\[
\frac{\delta_{0,\e}^{(N-2)/2}}{\sqrt\e}
=\sqrt{\mathtt d_N}d_{0,\e}
\longrightarrow\sqrt{\mathtt d_N}\bd d_0(\rho),
\qquad
\frac{\delta_{1,\e}^{(N-2)/2}}{\sqrt\e}
=\sqrt{\mathtt d_N}d_{1,\e}
\longrightarrow\sqrt{\mathtt d_N}\bd d_1(\rho),
\]
the limiting profile is precisely \eqref{eq:nodal-profile}.
\end{proof}

\begin{proposition}\label{prop:nodal-normal}
Let $N\geq13$ and $\rho\in\{\rho^*,\rho^{**}\}$.
Then the profile $\Phi_\rho$ defined in
\eqref{eq:nodal-profile} satisfies
\[
\partial_\nu\Phi_\rho(\theta)>0
\qquad\text{for every }\theta\in\partial\mathcal B.
\]
\end{proposition}

\begin{proof}
By \eqref{critrad} and \eqref{signo},
$\rho<1/\sqrt3$ and $\mathfrak t_3(\rho)=T(\rho)$.
Differentiating \eqref{Green} in the exterior normal direction,
we obtain
\[
\partial_{\nu_\theta}G(\theta,\xi)
=-(N-2)\frac{1-|\xi|^2}{|\theta-\xi|^N}.
\]
Hence, by \eqref{solsis},
\[
\partial_\nu\Phi_\rho(\theta)
=(N-2)\alpha_N\sqrt{\mathtt d_N}\,\bd d_1(\rho)
\left[
-T(\rho)+(1-\rho^2)
\sum_{j=1}^3\frac1{|\theta-\xi_j(\rho)|^N}
\right].
\]
Since $|\xi_j(\rho)|=\rho$ and
$\sum_{j=1}^3\xi_j(\rho)=0$, we have
\[
\frac13\sum_{j=1}^3|\theta-\xi_j(\rho)|^2=1+\rho^2.
\]
Jensen's inequality and Lemma~\ref{lem:T-bound} therefore give
\[
(1-\rho^2)\sum_{j=1}^3\frac1{|\theta-\xi_j(\rho)|^N}
\geq\frac{3(1-\rho^2)}{(1+\rho^2)^{N/2}}
>T(\rho),
\]
which proves the assertion.
\end{proof}

We can now conclude the proof of Theorem~\ref{thm:nodal}.

\begin{proof}[Proof of Theorem~\ref{thm:nodal}]
Let $i\in\{1,2\}$ and let $\rho\in\{\rho^*,\rho^{**}\}$ be the corresponding limiting radius. The profile $\Phi_\rho$ is of class $C^1$ in $\overline{\mathcal B}\setminus\{0,\xi_1(\rho),\xi_2(\rho),\xi_3(\rho)\}$. Hence, by Proposition~\ref{prop:nodal-normal} and the compactness of $\partial\mathcal B$, there exist $c>0$ and $\tau\in(0,1-\rho)$ such that
\[
\nabla\Phi_\rho(x)\cdot\frac{x}{|x|}\geq2c
\qquad\text{for every }1-\tau\leq|x|\leq1.
\]
By Lemma~\ref{lem:nodal-farfield} with $b=1-\tau$, there exists $\e_1>0$ such that, for every $\e\in(0,\e_1)$,
\[
\nabla u_\e^{(i)}(x)\cdot\frac{x}{|x|}\geq c\sqrt\e>0
\qquad\text{for every }1-\tau\leq|x|\leq1.
\]
In particular, $\partial_\nu u_\e^{(i)}>0$ on $\partial\mathcal B$. Moreover, since $u_\e^{(i)}=0$ on $\partial\mathcal B$, for every $x$ with $1-\tau<|x|<1$ we get
\[
u_\e^{(i)}(x)=-\int_{|x|}^1\nabla u_\e^{(i)}\Big(s\frac{x}{|x|}\Big)\cdot\frac{x}{|x|}\dd s\leq-c\sqrt\e\,(1-|x|)<0.
\]
Therefore the zeros of $u_\e^{(i)}$ in $\mathcal B$ are contained in $\{|x|\leq1-\tau\}$, which is a closed set, and the inclusion for $\mathcal Z(u_\e^{(i)})$ follows. Finally, it is enough to take the smaller of the two values of $\tau$ and of $\e_1$ obtained for $i=1,2$.
\end{proof}

Figure~\ref{fig:nodal} shows the sign of the limit profile $\Phi_\rho$ for $N=13$ at the two critical radii, which we computed numerically as the two zeros of $Q$. In both cases the nodal set is a closed hypersurface enclosing the positive bubble, and it stays far from $\partial\mathcal B$, as stated in Theorem~\ref{thm:nodal}.

\section{Nonexistence for \texorpdfstring{$k \geq 4$}{k>=4}}\label{nE4}

\begin{proof}[Proof of Theorem~\ref{thm:k-geq4}]
By \eqref{fkprimoutile},
\[
f_k'(\rho)
=\frac{k}{2}\bd d_1^2(\rho)
\bigl(A_k'(\rho)+2\beta'(\rho)\mathfrak t_k(\rho)\bigr).
\]

If $3\leq N\leq12$, Lemma~\ref{nonexgeq4}
gives $f_3'>0$ on $(\rho_{0,3},1)$.
For $\rho\in(\rho_{0,k},1)$,
Lemma~\ref{lem:polygon-estimates} gives
\[
A_3(\rho)\geq A_k(\rho)>0,\qquad
A_k'(\rho)\geq A_3'(\rho),\qquad
\mathfrak t_k(\rho)\leq\mathfrak t_3(\rho).
\]
Since $\beta'<0$, we obtain
\[
A_k'+2\beta'\mathfrak t_k
\geq A_3'+2\beta'\mathfrak t_3>0,
\]
and hence $f_k'>0$.

Let now $N\geq13$ and $\rho\in(\rho_{0,k},1)$. 
By Lemma~\ref{lem:polygon-estimates},
\[
\rho>\frac2{\sqrt{15}},
\qquad
A_k'>\frac{2(N-2)\rho}{1-\rho^2}A_k.
\]
Moreover, \eqref{eqtk} and $k\geq4$ imply
$3\beta\mathfrak t_k<A_k$.
Using $\beta'<0$ and
\[
\frac{\beta'}{\beta}
=-\frac{N-2}{\rho(1-\rho^{N-2})},
\]
we deduce
\[
A_k'+2\beta'\mathfrak t_k
>
2(N-2)A_k
\left[
\frac{\rho}{1-\rho^2}
-\frac1{3\rho(1-\rho^{N-2})}
\right].
\]
The expression in square brackets is positive, since
\[
\begin{aligned}
3\rho^2(1-\rho^{N-2})
&\geq3\rho^2(1-\rho^4)=3\rho^2(1+\rho^2)(1-\rho^2)>\frac{76}{75}(1-\rho^2)>1-\rho^2.
\end{aligned}
\]
Thus $f_k'>0$ also in this case.
Proposition~\ref{prop:link} then implies that $\Psi_k$
has no critical point in $\mathcal M$.
Hence, by Lemma~\ref{lem:limit-config}, there is no family of solutions as in the statement.
\end{proof}

%
\appendix
\section{Some technical lemmas}
The lemmas of this appendix are elementary, but we did not find them in the literature in the form needed here, so we give their proofs. Lemma~\ref{lem:positive-coefficients} is a roots-of-unity filter applied to the binomial series. The key step in Lemma~\ref{lem:auxiliary-monotonicity} is the classical fact that $\rho S'(\rho)/S(\rho)$ is nondecreasing for a power series $S$ with nonnegative coefficients, which follows from the Cauchy--Schwarz inequality. Lemmas~\ref{lem:T-bound} and \ref{lem:polygon-estimates} are explicit estimates on the functions $T$ and $A_k$.

\begin{lemma}\label{lem:positive-coefficients}
Let $N>2$. There exist coefficients $c_j\geq0$, $j\geq0$, such that
\[
(1-\rho^2)^{-(N-2)}+2(1+\rho^2+\rho^4)^{-\frac{N-2}{2}}=\sum_{j=0}^{\infty}c_j\rho^{2j},\qquad |\rho|<1,
\]
where the series converges absolutely, and $c_0=3$, $c_1=0$, $c_2=\frac{3(N-2)^2}{4}$.
\end{lemma}

\begin{proof}
Set $\gamma:=\frac{N-2}{2}$ and write $(1-z)^{-\gamma}=\sum_{m\geq0}b_mz^m$ for $|z|<1$, where $b_0=1$ and $b_m=\frac1{m!}\prod_{\ell=0}^{m-1}(\gamma+\ell)>0$ for $m\geq1$; in particular $b_1=\gamma$. Let $\omega:=e^{2\pi i/3}$. For real $|\rho|<1$ the numbers $1-\omega\rho^2$ and $1-\omega^{-1}\rho^2$ are complex conjugate, have positive real part, and their product is $1+\rho^2+\rho^4$; hence, using the principal branch of the power,
\[
(1-\rho^2)^{-2\gamma}+2(1+\rho^2+\rho^4)^{-\gamma}
=\sum_{h=0}^{2}(1-\omega^h\rho^2)^{-\gamma}(1-\omega^{-h}\rho^2)^{-\gamma}
=\sum_{m,m'\geq0}b_mb_{m'}\Big(\sum_{h=0}^{2}\omega^{h(m-m')}\Big)\rho^{2(m+m')}.
\]
Since $\sum_{h=0}^2\omega^{hq}$ equals $3$ if $3\mid q$ and $0$ otherwise, we get
\[
c_j=3\sum_{\substack{m+m'=j\\ 3\mid(m-m')}}b_mb_{m'}\geq0,\qquad
\sum_{j\geq0}c_j|\rho|^{2j}\leq3\Big(\sum_{m\geq0}b_m|\rho|^{2m}\Big)^2<\infty .
\]
The only admissible pairs $(m,m')$ are $(0,0)$ for $j=0$ and $(1,1)$ for $j=2$, and there is none for $j=1$. Thus $c_0=3$, $c_1=0$ and $c_2=3b_1^2=\frac{3(N-2)^2}{4}$.
\end{proof}

\begin{lemma}\label{lem:auxiliary-monotonicity}
Let $a,p,R>0$ and let $S(\rho)=\sum_{j\geq0}a_j\rho^{p+2j}$ for $0<\rho<R$, where $a_j\geq0$, at least one $a_j$ is positive, and the series converges on $(0,R)$. Set $T:=a+S$ and
\[
\Theta(\rho):=1-\rho^p\big(1-T(\rho)\big)-2p\,\frac{T(\rho)}{\rho T'(\rho)} .
\]
Then $\Theta'(\rho)>0$ for every $\rho\in(0,R)$ such that $\rho^p<8a$.
\end{lemma}

\begin{proof}
The series can be differentiated term by term, and $S,S'>0$ on $(0,R)$. Writing $\lambda_j:=p+2j$, the Cauchy--Schwarz inequality gives
\[
(\rho S')^2=\Big(\sum_{j}\lambda_ja_j\rho^{\lambda_j}\Big)^2\leq\Big(\sum_{j}a_j\rho^{\lambda_j}\Big)\Big(\sum_{j}\lambda_j^2a_j\rho^{\lambda_j}\Big)=S\,\rho(\rho S')' .
\]
Since $T'=S'$ and $T=a+S$, it follows that
\[
\rho\Big(\frac{T}{\rho T'}\Big)'=1-\frac{T\,\rho(\rho S')'}{(\rho S')^2}\leq1-\frac TS=-\frac aS .
\]
Using also $S(1-T)=S(1-a-S)\leq S(1-S)\leq\frac14$, we obtain
\[
\rho\Theta'=-p\rho^p(1-T)+\rho^{p+1}T'-2p\rho\Big(\frac{T}{\rho T'}\Big)'
\geq\frac pS\big(2a-\rho^pS(1-T)\big)+\rho^{p+1}T'
\geq\frac pS\Big(2a-\frac{\rho^p}{4}\Big)+\rho^{p+1}T',
\]
which is positive whenever $\rho^p<8a$.
\end{proof}

\begin{lemma}\label{lem:T-bound}
Let $N\geq3$ and let $T$ be defined by \eqref{T}. Then
\[
T(\rho)<M(\rho):=\frac{3(1-\rho^2)}{(1+\rho^2)^{N/2}}\qquad\text{for every }0<\rho\leq\tfrac1{\sqrt3}.
\]
\end{lemma}

\begin{proof}
By \eqref{T}, $T(\rho)<\rho^N(1-\rho^2)^{1-N}+3^{-(N-2)/2}$. Since $\rho^2\leq\frac13$, we have $\rho^2(1+\rho^2)\leq(1-\rho^2)^2$ and $1+\rho^2\leq\frac43$, so that
\[
\frac{\rho^N}{(1-\rho^2)^{N-1}M(\rho)}=\frac13\Big(\frac{\rho\sqrt{1+\rho^2}}{1-\rho^2}\Big)^N\leq\frac13,\qquad
\frac{3^{-(N-2)/2}}{M(\rho)}=\frac1{1-\rho^2}\Big(\frac{1+\rho^2}{3}\Big)^{N/2}\leq\frac32\Big(\frac49\Big)^{N/2}=\Big(\frac23\Big)^{N-1}.
\]
Therefore $T/M<\frac13+(\frac23)^{N-1}\leq\frac79<1$ for $N\geq3$.
\end{proof}

\begin{lemma}\label{lem:polygon-estimates}
Let $N\geq3$ and $k\geq3$. For every $\rho\in(0,1)$, $A_k(\rho)\leq A_3(\rho)$ and $A_k'(\rho)\geq A_3'(\rho)$. If $A_k(\rho)>0$, then
\[
\mathfrak t_k(\rho)\leq\mathfrak t_3(\rho)\qquad\text{and}\qquad A_k'(\rho)>\frac{2(N-2)\rho}{1-\rho^2}A_k(\rho).
\]
Furthermore, if $N\geq4$ and $k\geq4$, then $\rho_{0,k}>\frac2{\sqrt{15}}$.
\end{lemma}

\begin{proof}
For $0<\rho<1$ and $0<s\leq2$ set $D(\rho,s):=(1-\rho^2)^2+\rho^2s^2$ and
\[
g(\rho,s):=(\rho s)^{2-N}-D(\rho,s)^{-\frac{N-2}{2}},\qquad\text{so that, by \eqref{Akk},}\qquad
A_k(\rho)=(1-\rho^2)^{2-N}-\sum_{j=1}^{k-1}g\Big(\rho,2\sin\frac{\pi j}{k}\Big).
\]
Clearly $g>0$, and a direct differentiation gives
\[
\partial_sg=-(N-2)\rho^2s\Big[(\rho s)^{-N}-D^{-\frac N2}\Big]<0,\qquad
-\partial_\rho g=\frac{N-2}{\rho}\Big[g+(1-\rho^4)D^{-\frac N2}\Big]>0 .
\]
Thus both $g$ and $-\partial_\rho g$ are positive and strictly decreasing in $s$. Since $2\sin\frac\pi k=2\sin\frac{\pi(k-1)}{k}\leq\sqrt3$, keeping only the terms $j=1$ and $j=k-1$ in the sum above and in its $\rho$-derivative we get $A_k\leq A_3$ and $A_k'\geq A_3'$.

If $A_k(\rho)>0$, then $A_3(\rho)>0$, and by \eqref{eqtk} $\mathfrak t_k$ is the positive root of $t^2+(k-1)\beta t=A_k$; this root is increasing in $A_k$ and decreasing in $k$, so $\mathfrak t_k\leq\mathfrak t_3$. Moreover, the positivity of $g$ and $-\partial_\rho g$ gives $A_k<(1-\rho^2)^{2-N}$ and $A_k'>2(N-2)\rho(1-\rho^2)^{1-N}$, whence $A_k'>\frac{2(N-2)\rho}{1-\rho^2}A_k$.

Finally, let $N\geq4$, $k\geq4$ and $r_c:=2/\sqrt{15}$. Since $2\sin\frac\pi k\leq\sqrt2$ and the sum contains at least one more positive term, $A_k(\rho)<(1-\rho^2)^{2-N}-2g(\rho,\sqrt2)$. Using $(1-r_c^2)^2=\frac{121}{225}$, $2r_c^2=\frac{120}{225}$ and the fact that $\gamma\mapsto x^\gamma-y^\gamma$ is increasing for $x>1>y>0$, we obtain, with $\gamma=\frac{N-2}{2}\geq1$,
\[
(1-r_c^2)^{N-2}A_k(r_c)<1-2\Big[\Big(\frac{121}{120}\Big)^{\gamma}-\Big(\frac{121}{241}\Big)^{\gamma}\Big]\leq1-2\Big(\frac{121}{120}-\frac{121}{241}\Big)<0 .
\]
Since $A_k$ is strictly increasing and $A_k(\rho_{0,k})=0$, it follows that $\rho_{0,k}>r_c$.
\end{proof}

\section*{Data Availability Statements}
All data generated or analysed during this study are included in this article.
\section*{Declarations}
\textbf{Conflicts of Interest:} The authors declare they have no financial interests.
\section*{Acknowledgements}
During the preparation of this work, the authors used ChatGPT (OpenAI) and Claude (Anthropic) to assist with selected calculations and to revise the manuscript. The numerical computations and the plots in Figure~\ref{fig:nodal} were produced with Python code generated with the help of ChatGPT; the code is available from the authors upon request. All mathematical statements, proofs, and conclusions were rigorously established by the authors, who take full responsibility for the accuracy and originality of the entire work.

\end{document}